\documentclass{article}

\usepackage{graphicx}%
\usepackage{multirow}%
\usepackage{amsmath,amssymb,amsfonts,bm}%
\usepackage{amsthm}%
\usepackage{mathrsfs}%
\usepackage[title]{appendix}%
\usepackage{caption}
\usepackage{xcolor}%
\usepackage{textcomp}%
\usepackage{manyfoot}%
\usepackage{booktabs}%
\usepackage{algorithm}%
\usepackage{placeins}
\usepackage{algorithmicx}%
\usepackage{algpseudocode}%
\usepackage{listings}%

\usepackage[colorlinks=true]{hyperref}
\usepackage{pdfsync}

\newtheorem{remark}{Remark}
\newtheorem{definition}{Definition}%

\title{Numerical solving of an optimal control problem in large time horizon: the aerial vehicle guidance}

\author{
Veljko A\v{s}kovi\'{c}\footnote{MBDA France, 92350 Le Plessis Robinson, France (\texttt{veljkoaskovic@hotmail.com}).}
\and
Emmanuel Tr\'elat\footnote{Sorbonne Universit\'e, Universit\'e Paris Cit\'e, CNRS, Inria, Laboratoire Jacques-Louis Lions, LJLL, F-75005 Paris, France (\texttt{emmanuel.trelat@sorbonne-universite.fr}).}
\and
Hasnaa Zidani\footnote{LMI, Insa Rouen Normandie, 76800 Rouen, France (\texttt{hasnaa.zidani@insa-rouen.fr}).}
}

\date{}

\begin{document}
\maketitle


\abstract{In this paper we consider an optimal control problem in large time horizon and solve it numerically. More precisely, we are interested in an aerial vehicle guidance problem: launched from a ground platform, the vehicle aims at reaching a ground/sea target under specified terminal conditions while minimizing a cost modelling some performance and constraint criteria. Our goal is to implement the indirect method based on the Pontryagin maximum principle (PMP) in order to solve such a problem. After modeling the problem, we implement continuations in order to ``connect'' a simple problem to the original one. Particularly, we exploit the turnpike property in order to enhance the efficiency of the shooting.}

\bigskip




\maketitle

\section{Introduction}\label{sec1}


The study of autonomous air-breathing vehicles is of significant interest in both military and civilian applications, particularly when considering constraints related to  maneuverability and operational efficiency. In this paper, we focus on the trajectory optimization of a self-guided, pilotless vehicle. The vehicle is designed to follow a pre-programmed path. The primary objective is to ensure that the vehicle reaches a designated target while adhering to various operational constraints.

A key challenge in trajectory planning is to minimize detectability by enemy radar systems, which is typically achieved by maintaining a low-altitude cruise phase. Additionally, the vehicle must arrive at the target with a predefined flight path angle and, potentially, a specified terminal velocity. Given these constraints, it is natural to expect a structured trajectory composed of three distinct phases: an initial transitory phase, a cruise phase, and a bunt phase leading to the impact on the target.

The main contribution of this work is the formulation and numerical resolution of an optimal control problem governing the vehicle's trajectory in the vertical plane. This problem incorporates physical and operational constraints while ensuring adherence to expected qualitative behaviors. A particular emphasis is placed on the application of the Shooting Method, a powerful numerical approach for solving boundary value problems. However, its implementation in this context presents specific challenges, including the need to maintain stealth during the cruise phase and the inherent non-smooth nature of the vehicle's dynamics.

To address these difficulties, we propose strategies to enhance the initialization of the Shooting Method and techniques to stabilize the computations. The effectiveness of the approach is demonstrated through numerical simulations, which highlight the performance of the Shooting Method and the advantages of the proposed refinements in optimizing the vehicle's trajectory.

\section{Optimization problem}\label{sec2}
\subsection{Context of the study}
We consider in this paper a self guided, pilotless, continuously powered air-breathing vehicle that flies like a miniature airplane, supported by its aerodynamic surfaces. Launched from a dedicated platform and propelled by a turbojet, it aims at reaching a specified location on the ground (defined through its coordinates).

\noindent Prior to the flight, the desired trajectory is pre-programmed on the onboard computer. Then it is followed by the vehicle during the flight updating its position thanks to the navigation sensors (Global Positioning System \& Inertial Nagivation Unit). If necessary, the actual trajectory can be adjusted/corrected with respect to the specified one.
\vspace{0.1cm}

\noindent Moreover, the specified trajectory has to respect some performance/constraint criteria:
\begin{itemize}
\item The vehicle has to be as discrete as possible with respect to ennemy sensors (radars). To do so, a classical strategy consists in flying at a low altitude, that will be named ``cruise altitude'' in the sequel of the paper. The cruise altitude, denoted by $h_{c}$, is usually comprised between $100$m and $300$m with respect to the Earth surface.
\item The vehicle has to reach the target with a specified flight path angle, denoted by $\gamma_{f}<0$ (diving interception) and possibly a specified terminal speed $v_{f}>0$. 
\item Due to the stealth constraint, it is desirable to accomplish the mission in a minimum of time. 
\end{itemize}

\noindent In light of the preceding discussion, it is reasonable to expect that the trajectory in the vertical plane would exhibit a behavior similar to the one illustrated below.
\begin{figure}[h]
\begin{center}
\captionsetup{justification=centering}
\includegraphics[width=0.8\linewidth]{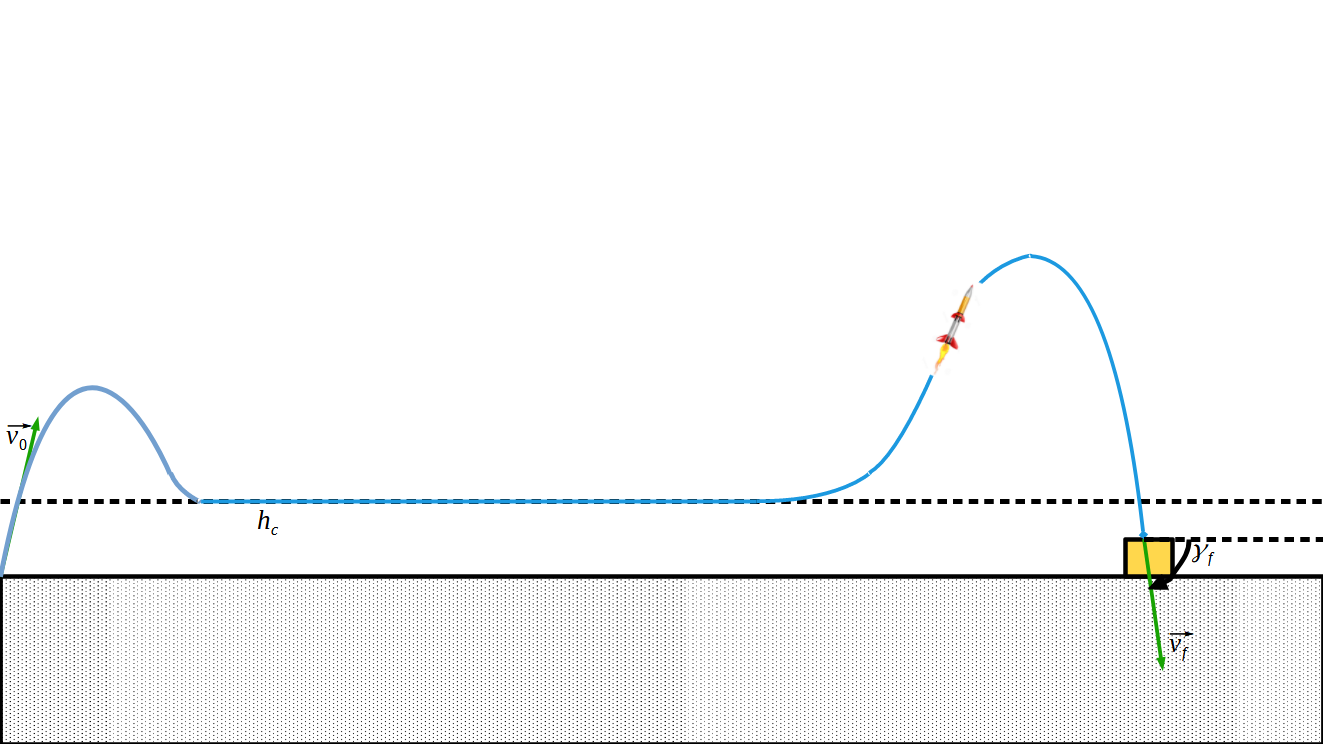}
\caption{Expected trajectory in the vertical plane}
\label{fig1}	
\end{center}
\end{figure}
\FloatBarrier
\noindent  Consequently, the trajectory can likely be divided into three segments:
\begin{enumerate}
\item \emph{Initial transitory phase}: The vehicle transitions from its launching position to the cruise altitude.
\item \emph{Cruise phase}: The vehicle maintains an altitude close to $h_{c}$.  Actually, given the visibility constraints, one can intuitively expect the vehicle to remain "stuck" to the cruise altitude $h_{c}$ most of its flight time. This behavior ensures that the vehicle remains within an optimal range for observation or engagement while adhering to operational constraints.
\item \emph{Bunt phase}: the vehicle leaves the cruise altitude et climbs up to the peak of the trajectory before "diving" onto the target. The maximum altitude reached during this phase is expected to depend on the desired terminal flight path angle and potentially on the prescribed terminal velocity.
\end{enumerate}
\noindent Depending of the initial (resp. terminal) speeds and flight path angles,  the first and final phases can possibly exhibit a "parabolic shape" as illustrated on the figure \ref{fig1}.
\vspace{0.1cm}

\noindent In this paper, our primary focus is the computation of such a trajectory in the vertical plane, ensuring that it satisfies the imposed physical and operational constraints while adhering to the expected qualitative behavior, as well as the formulation of the optimal control problem that generates it, in order to gain a deeper understanding of the constraints and the performance optimized by this type of trajectory.
\vspace{0.1cm}

\noindent Once the optimal control problem is formulated, we will focus on numerical methods for its resolution. In particular, we will consider the implementation of the Shooting Method. This approach, known for its high efficiency and accuracy, generally requires good initialization to ensure convergence to an optimal solution.
Its application to the model considered here, however, presents specific challenges, including the non-visibility requirement (during the cruise phase), the non-smooth nature of the dynamics. These difficulties will be discussed in detail in the following sections.
We will also present the proposed solutions to overcome these challenges, including strategies to improve initialization and techniques to stabilize the computations. Finally, the results will be illustrated through numerical simulations, highlighting the performance of the Shooting Method in this specific context and the advantages brought by the suggested adaptations.

\subsection{Turnpike property}
\noindent Seen from afar, the expected trajectory will remain most of the time ``close to the cruise altitude'' which, from the optimal control problem formalism, reminds of  a so-called \emph{turnpike property}: roughly speaking, it reflects the fact that, for suitable optimal control problems in a sufficiently large time horizon, any optimal solution thereof remains, most of the time, close to the optimal solution of an associated static optimization problem. This optimal static solution is referred to as the turnpike (the name stems from the idea that a turnpike is the fastest route between two points which are far apart, even if it is not the most direct route).
\vspace{0.2cm}

\noindent The turnpike phenomenon was first observed and investigated in the context of discrete-time optimal control problems by economists (see \cite{Neu1945}, \cite{DorHut1958}). Over the years, the concept has been rigorously formalized, leading to several variants, some of which are stronger than others (see \cite{Zas2015}). For instance, the \emph{exponential turnpike property}, established in \cite{TreZua2015}, ensures that the entire optimal triple—comprising the state, adjoint, and control—remains exponentially close to the optimal solution of the corresponding static controlled problem, except near the initial and terminal times of the horizon.
\vspace{0.2cm}

\noindent From the numerical point of view, the turnpike notion is of significant practical relevance, offering valuable insights into the design and computational reconstruction of optimal controllers. For example, the turnpike property has been exploited to ensure the stability of numerical optimization algorithms, particularly in the context of model predictive control schemes \cite{Gru2022}. Furthermore, \cite{TreZua2015} and \cite{CaiFerTreZid2022} illustrate how the turnpike phenomenon can be exploited in order to improve the initialization of shooting methods, which are highly sensitive to initial guesses, especially in long-horizon problems.
\vspace{0.2cm}

\noindent For a comprehensive survey on the turnpike property and its applications in optimal control problems, see, for example, \cite{Gru2022bis}. This body of work highlights the theoretical significance and the computational benefits of understanding and exploiting the turnpike property in various settings.

\section{Equations of the motion}\label{sec3}
\subsection{Assumptions}
In what follows, we make the ``flat earth''assumption (the traveled distance being negligible compared to the radius of Earth). Moreover we assume that the launcher and target are located at sea level altitude thus one sets, for sake of simplicity: $(h_{0}, h_{f}):=(0,0)$.
\subsection{Coordinate systems}\label{subsec3}
\begin{figure}[h]
\begin{center}
\captionsetup{justification=centering}
\centerline{\includegraphics[width=0.8 \linewidth]{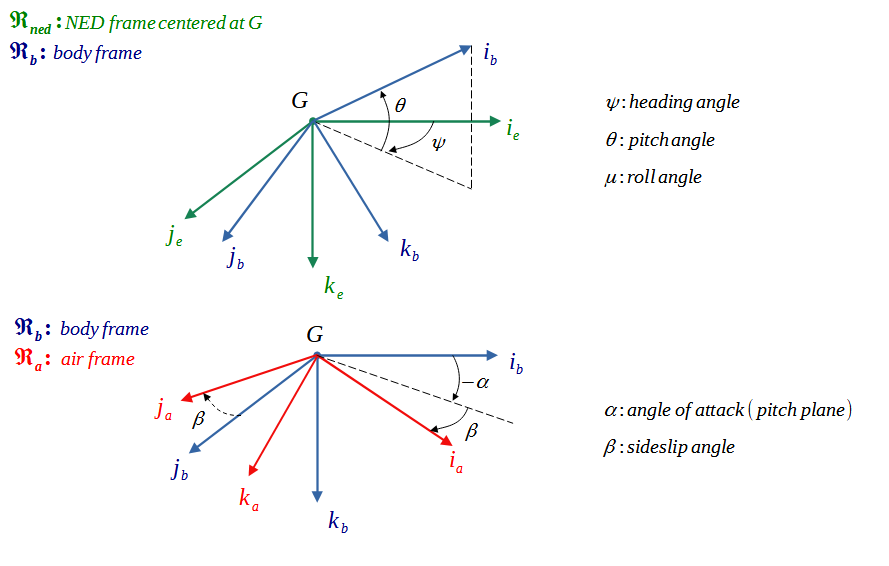}}
\caption{Coordinate systems}
\label{fig2}	
\end{center}
\end{figure}

The derivation of the equations of motion is clarified by defining a number of coordinate systems. The three coordinate systems used here are the following, all of them being direct and orthonormal (see figure \ref{fig2}):
\begin{enumerate}
\item The local earth frame $\mathcal{R}_{ned}=(E,\textbf{i}_{e},\textbf{j}_{e},\textbf{k}_{e})$ is fixed to the surface of the  earth at mean sea level: $\textbf{i}_{e}$ is oriented towards the local North, $\textbf{j}_{e}$ towards the East and $\textbf{k}_{e}$ is down oriented . This referential is usually called "North-East Down" (NED).
\item The body frame $\mathcal{R}_{b}=(G, \textbf{i}_{b},\textbf{j}_{b},\textbf{k}_{b})$ is conventional to the body of the vehicle. The center of this frame is at the center of gravity of the vehicle, and its components are forward, out of the right side, and down. The orientation is calculated from the $\mathcal{R}_{ned}$ frame through three successive rotations: $\psi$-heading angle around $z$ axis, $\theta$-pitch angle around $y$ axis and $\mu$-roll angle around $x$ axis.  
\item The air frame $\mathcal{R}_{a}=(G, \textbf{i}_{a},\textbf{j}_{a},\textbf{k}_{a})$ centered at the center of gravity and where $\textbf{i}_{a}$ is coincident with the velocity vector $\textbf{v}$.  The orientation is obtained from the body frame through two successive rotations: $\alpha$-angle of attack around $y$ axis, $\beta$-sideslip angle around $z$ axis. 
\end{enumerate}

\FloatBarrier

\begin{remark} \label{rem1} As one wants the angle of attack to be positive when the body is above the velocity in the pitch plane, we change the sign convention accordingly  (see the figure \ref{fig2}). 
\end{remark}

\subsection{Kinematic equations}\label{subsec4}
Kinematics is used to derive the differential equations for the coordinates of G, denoted by $(x,y,-h)$ $h$ being the altitude of the vehicle above the sea level. The basic relation is:
\begin{equation}
\dfrac{d\textbf{EG}}{dt}=\dot{x}\textbf{i}_{e}+\dot{y}\textbf{j}_{e}-\dot{h}\textbf{k}_{e}=\textbf{V}=v\textbf{i}_{a}
\label{kin}
\end{equation}
The coordinates of the vector $\textbf{X}$ in the air frame $\mathcal{R}_{a}$ can be expressed in the NED frame $\mathcal{R}_{ned}$ as follows:
\begin{equation*} 
\textbf{X}|_{\mathcal{R}_{ned}}=P_{\mathcal{R}_{ned}/\mathcal{R}_{a}}\textbf{X}|_{\mathcal{R}_{a}}=P_{3}(\psi)P_{2}(\theta)P_{1}(\mu)P_{2}(-\alpha)P_{3}(\beta)\textbf{X}|_{\mathcal{R}_{a}}
\end{equation*}
where:
\begin{align*}
&P_{1}(\mu)=\begin{pmatrix}
1 & 0 & 0 \\
0 & \cos \mu & \sin \mu \\
0 & -\sin \mu & \cos \mu
\end{pmatrix}, P_{2}(\theta)=\begin{pmatrix}
\cos \theta & 0 & \sin \theta \\
0 & 1 & 0 \\
-\sin \theta & 0 & \cos \theta
\end{pmatrix},  P_{3}(\psi)=\begin{pmatrix}
\cos \psi & -\sin\psi & 0 \\
\sin \psi & \cos \psi & 0 \\
0 & 0 & 1
\end{pmatrix}
\end{align*}
With our notations, \eqref{kin} reads:
\begin{equation}
\left(\begin{array}{ccc}
\dot{x} \\
\dot{y} \\
-\dot{h}
\end{array}
\right)=P_{3}(\psi)P_{2}(\theta)P_{1}(\mu)P_{2}(-\alpha)P_{3}(\beta)\left(\begin{array}{ccc}
v \\
0 \\
0
\end{array}
\right)
\label{kin1}
\end{equation}

\subsection{Dynamic equations}\label{subsec5}
Newton’s second law can be expressed as follows:
\begin{equation}
m\dot{\textbf{V}}+\dot{m}\textbf{V}=\textbf{T}+\textbf{L}+\textbf{D}+\textbf{W},
\label{eqdyn}
\end{equation}
where $\textbf{V}$ is the velocity, $m$ is the mass of the vehicle, $\textbf{T}$ is the thrust, $\textbf{D}$ is the aerodynamic drag, $\textbf{L}$ is the aerodynamic lift, and $\textbf{W}$ is the weight.
\vspace{0.2cm}
	
\noindent \textbf{Weight}: The weight is given by $\textbf{W}=mg\textbf{k}_{e}$, where $g$ is the gravitational acceleration and $\textbf{k}_{e}$ is the unit vector pointing downward in the local Earth-fixed frame.
 Moreover, by introducing the specific fuel consumption denoted by $C_{s}$, the mass equation reads:
\begin{equation}
\dot{m}=-C_{s}T
\label{meq}
\end{equation}
where $T$ denotes the magnitude of the thrust force.
\vspace{0.1cm}

\noindent \textbf{Drag and Lift}: The drag force acting on the vehicle is the frictional force exerted by the air, directed opposite to the velocity vector $\textbf{V}$. Typically, the effect of the drag force is summarized using a single coefficient, $C_{D}$, called the drag coefficient. This coefficient depends on several factors, including the angle of attack (the angle between the velocity vector $\textbf{V}$ and the principal body axis $\textbf{i}{b}$), the Mach number, and the aerodynamic configuration of the vehicle.
 The drag force is mathematically expressed as:
\begin{equation}
\textbf{D}=-D\textbf{ i}_{a} \text{ with } D=\bar{q}(h,v)SC_{D}(h,v,\alpha,\beta)
\label{drag}
\end{equation}
with $S=\dfrac{\pi d^{2}}{4}$ being the vehicle's reference cross-sectional area, $\bar{q}(h,v)=\frac{1}{2}\rho(h).v^{2}$ is the dynamic pressure. Here, $\rho(h)$ denotes the air density, $v$ is the velocity magnitude, and $h$ is the altitude.
\vspace{0.1cm}

\noindent For altitudes below $20$ km, the air density can be analytically approximated by the formula:
\begin{equation}
\rho(h)=\rho_{0}\exp(-h/h_{r})
\label{rhofor}
\end{equation}
where $\rho_{0}$ is the air density of the standard atmosphere at the sea level, and $h_{r}$ is a fixed reference altitude.This exponential model provides a convenient approximation of the atmosphere and will be used throughout this study. For further details on the atmospheric model, we refer the reader to \cite{Siouris2003}.
\vspace{0.2cm}
22
\noindent The lift is defined as the aerodynamic force that acts orthogonally to the velocity vector. The analytical expression of the lift is:
\begin{equation}
\textbf{L}=L\textbf{i}_{a} \wedge (\textbf{i}_{b} \wedge \textbf{i}_{a}) \text{ with } L=\bar{q}(h,v)SC_{L}(h,v,\alpha,\beta)
\label{lift}
\end{equation}
The scalar coefficient $C_{L}$ called lift coefficient has the same dependencies than the drag coefficient $C_{D}$.
\vspace{0.1cm}

\noindent \textbf{Thrust}: In terms of direction, the thrust vector is aligned with the vehicle's principal body axis and directed forward, thus $\textbf{T}=T\bf{i}_{b}$. It is also known that thrust and specific fuel consumption $C_{s}$ satisfy functional relations of the form:
\begin{equation}
T = T(h, v, u_{1}), \hspace{0.5cm} C_{s} = C_{s}(h, v, u_{1})
\label{thrust}
\end{equation}
where $u_{1}$ is the thrust throttle coefficient.
\vspace{0.1cm}

\noindent For more details on the forces exerted on the vehicle, see \cite{Siouris2003}, \cite{Fleeman2012} and \cite{Zarchan2012}.
\subsection{Flight in the vertical plane} \label{subsec6}
Within the frame of our study, we restrict ourselves to the flight of the vehicle in the vertical plane. Consequently, from now on, the heading, the roll and the sideslip angles are such that $\psi=0$, $\mu=0$ and $\beta=0$.

\begin{figure}[h!]
\begin{center}
\captionsetup{justification=centering}
\includegraphics[width=0.7\linewidth]{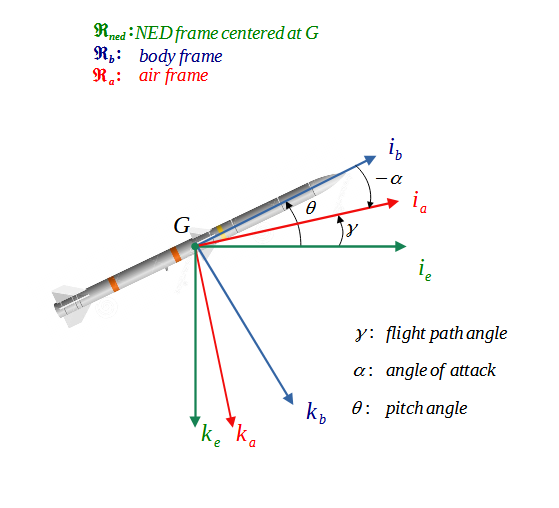}
\caption{Forces acting on the vehicle in flight}
\label{fig3}	
\end{center}
\end{figure}
\FloatBarrier

\noindent Under the previous assumptions, introducing the flight path angle $\gamma:=\theta-\alpha$, \eqref{kin1} reduces to:
\begin{subequations} \label{kin2}
\begin{align}
&\dot{x}=v\cos \gamma  \label{kin21} \\ 
&\dot{h}=v\sin \gamma  \label{kin22}
\end{align}
\end{subequations}

\begin{remark} \label{rem2} The flight path angle $\gamma$ is classically the angle between the velocity vector and the horizontal reference in the vertical plane.
\end{remark}
\vspace{0.1cm}
\noindent By projecting \eqref{eqdyn} respectively on $\textbf{i}_{a}$ and $\textbf{i}_{b}$, one obtains:
\begin{subequations} \label{eqdyn1}
\begin{align}
&\dot{v}=\dfrac{T\cos\alpha-D-\dot{m}v}{m}-g\sin\gamma  \label{eqdyn11} \\
&\dot{\gamma}=\dfrac{L+T\sin\alpha}{mv}-\dfrac{g\cos\gamma}{v} \label{eqdyn12}
\end{align}
\end{subequations}

\subsection{Discussion of the model} \label{subsec7}
Gathering \eqref{meq}, \eqref{kin2} and \eqref{eqdyn1} one obtains:
\begin{subequations} \label{eqkindyn}
\begin{align}
&\dot{x}=v\cos \gamma  \label{eqkindyn1} \\
&\dot{h}=v\sin \gamma \label{eqkindyn2} \\
&\dot{v}=\dfrac{T(h,v,u_{1})\cos\alpha-D(h,v,\alpha)-\dot{m}v}{m}-g\sin\gamma  \label{eqkindyn3} \\
&\dot{\gamma}=\dfrac{L(h,v,\alpha)+T(h,v,u_{1})\sin\alpha}{mv}-\dfrac{g\cos\gamma}{v} \label{eqkindyn4} \\
& \dot{m}=-C_{s}T \label{eqkindyn5}
\end{align}
\end{subequations}
Consequently, the state variables are $x$, $h$, $v$, $\gamma$ and $m$ whereas the thrust throttle $u_{1}$ and the angle of attack $\alpha$ appear naturally as the control variables.
\vspace{0.1cm}

\noindent Considering endo-atmospheric applications, the vehicle can be stabilized only if the angle of attack $\alpha$ remains lower than a critical threshold $\alpha_{\max}$, which is assumed to be quite low for such type of vehicle (typically $\alpha_{\max} \approx 15^{\circ}$).  For this reason we will assume that in \eqref{eqkindyn} one has $\cos \alpha \approx 1$ and $\sin \alpha \approx 0$. 

\noindent Moreover, in the case where $\alpha \ll 1$, a Taylor expansion of the aerodynamic coefficients allows to obtain the following parabolic drag polar approximation:
\begin{equation}
 C_{D}= C_{D_{0}}+k_{c}.C_{L}^2
\label{paradragpolar}
\end{equation}
where $C_{D_{0}}$ is the zero-lift drag coefficient, $k_{c}.C_{L}^2$ is the induced drag coefficient and $k_{c}$ is the induced-by-the-lift drag factor. Theoretically, $C_{D_{0}}$ and $k$ depend of the Mach number. However, as the Mach number range covered by the vehicle is reduced, we assume that the latter coefficients are constant.
\vspace{0.1cm}

\noindent With the drag polar approximation \eqref{paradragpolar}, one remarks that the aerodynamic forces $D$ and $L$ can now be expressed as functions of $C_{L}$, as the angle of attack $\alpha$ ``disappears''. Consequently, the lift coefficient $C_{L}$ becomes naturally the control in place of the angle of attack $\alpha$.
\vspace{0.1cm}

\noindent Related to the thrust, the dependence on the parameters $h$, $v$ can be quite intricate: indeed it is function of the altitude, Mach number, surface and geometry of the air inlet entrance and outlet exit, the compressor pressure ratio... An overview of the existing technology and some empirical formulas to estimate the thrust performance are available for instance in \cite{Fleeman2012}. In our study, as we do not expect meaningful  variations in altitude et Mach number, we consider the following (somehow idealized) formula:
\begin{equation}
T  \approx T_{\max}u_{1}
\label{thrust1}
\end{equation}
where the throttle coefficient $u_{1} \in [\eta_{\min}, 1]$, and $T_{\max}$ is the maximum available thrust that can be delivered by the engine.
\vspace{0.1cm}

\noindent Finally one considers the following vehicle dynamical model:

 \begin{subequations} \label{eqkindynapp}
\begin{align}
&\dot{x}=v\cos \gamma  \label{eqkindynapp1} \\
&\dot{h}=v\sin \gamma \label{eqkindynapp2} \\
&\dot{v}=\dfrac{T_{\max}(1+C_{s}v)u_1-D(h,v,C_{L})}{m}-g\sin\gamma  \label{eqkindynapp3} \\
&\dot{\gamma}=\dfrac{L(h,v,C_{L})}{mv}-\dfrac{g\cos\gamma}{v} \label{eqkindynapp4} \\
& \dot{m}=-C_{s}T_{\max}u_1 \label{eqkindynapp5}
\end{align}
\end{subequations}

\noindent In this problem, the state variable is five-dimensional, while the control input consists of two possible components: the thrust throttle, which regulates the velocity and the lift coefficient $C_{L}$, which shapes the trajectory.
\vspace{0.2cm}

 \noindent \textbf{Assumptions}: In the rest of the paper, we slightly simplify the original problem by assuming:
\begin{itemize}
\item the thrust remains at it maximum value, i.e, $\forall t \in [0, t_{f}]$, $u_{1}(t)=1$
\item the induced-by-the-lift drag factor $k_{cz}$ is negligible thus $k_{cz} \approx 0$. Indeed, in practice, the coupling between the lift and drag force can be neglected at first order. 
\end{itemize}
\begin{remark} \label{rem4}
Under these assumptions, the differential equation governing the state variable $v(\cdot)$ no longer explicitly depends on $u_{2}$. 
\end{remark}
\vspace{0.2cm}

\noindent Finally, the lift coefficient $u_{2}$ becomes the sole control variable of our system. For simplicity of notations, it will be referred to as $u$ throughout the rest of the paper.
\vspace{0.2cm}

\noindent The control input is bounded due to the vehicle aerodynamic limitations: $|u(t)| \leqslant u^{m}$. In what follows, we denote $\xi:=\left(x \hspace{0.2cm} h \hspace{0.2cm} v \hspace{0.2cm} \gamma \hspace{0.2cm} m \right)^{T}$ and $U:= [-u^{m}, u^{m}]$.
\section{Optimal control problem setting}\label{sec4}
\subsection{Constraints and cost modelling} \label{subsec8}
 The constraints and performance described in Section \ref{sec2} must be integrated into the model through an appropriate objective function. A common (but not exclusive!) approach to account for the cruise altitude constraint is to include it as a penalization term in the cost function. Let us consider a Lagrangian cost composed of two terms: the first term represents the time of flight, while the second penalizes deviations from the cruise altitude $h_{c}$. The optimal control problem can then be formulated as:
\begin{equation}
(\textbf{OCP})_{k_{0},k_{1}}\left\{
\begin{aligned}
&\underset{(t_{f}>0, \hspace{0.02cm} u \in \mathcal{U})}{\min} J_{k_{0}, k_{1}}(t_{f}, u):=\displaystyle{\int_{0}^{t_{f}}} f^{0}(\xi(s), u(s))ds \\
&\dot{\xi}(s)=f(\xi(s),u(s)) \hspace{1cm} \forall s \in [0, t_{f}] \\
&u(s) \in U \hspace{1cm} \forall s \in [0, t_{f}] \\
& \xi(0)=\xi_{0}, \hspace{0.2cm} \xi(t_{f})=\xi_{f}
\end{aligned}
\right.
\label{contopt1}
\end{equation}
where:
\begin{itemize}
\item $\xi=\left(x \hspace{0.2cm} h \hspace{0.2cm}  v \hspace{0.2cm} \gamma \hspace{0.2cm} m 	\right)^{T}$ is the state
\item $f(\xi,u):=\left(\begin{array}{ccccc}
		v\cos \gamma \\
		v\sin \gamma \\
		\dfrac{T_{\max}(1+C_{s}v)u_{1}-D(h,v,u_{2})}{m}-g\sin\gamma \\
		\dfrac{L(h,v,u_{2})}{mv}-\dfrac{g\cos\gamma}{v} \\
		-C_{s}T_{\max}u_{1}
	\end{array}
	\right)
	$ is the dynamics
\item $f^{0}(\xi, u):=k_{0}+k_{1}\dfrac{(h(s)-h_{c})^2}{h_{c}^2}$ is the running cost
\item $\mathcal{U}$ is a set of admissible control strategies defined as $\mathcal{U}:=\{u: [0, +\infty[ \rightarrow U\}$ with $U$ defined above
\item $\xi_{0}=\left(x_{0} \hspace{0.2cm} h_{0} \hspace{0.2cm}  v_{0} \hspace{0.2cm} \gamma_{0} \hspace{0.2cm} m_{0}
	\right)^{T}$ and $\xi_{f}=\left(x_{f} \hspace{0.2cm} h_{f} \hspace{0.2cm}  v_{f} \text{ or } *\hspace{0.2cm} \gamma_{f} \text{ or }*  \right)^{T}$ are the prescribed initial and final states. $*$ designates the case where the state component is let free.
\item $(k_{0}, k_{1}) \in \mathbb{R}_{+}^{2}$ is the weight couple in the performance index. $k_{0}$ and $k_{1}$ weight respectively the time of flight and the penalization of the cruise altitude.
\end{itemize}
We have now a family of optimal control problems parametrized by the weights $k_{0}$ and $k_{1}$.
\subsection{Comments} \label{subsec9}
If one chooses $k_{0}>0$ and $k_{1}=0.0$,  then $(\textbf{OCP})_{k_{0},0}$ is the time minimum problem. If  $k_{1}>0$ and is ``large enough'', then the cruise altitude constraint becomes active and the optimal trajectory should correspond rather to the one illustrated in the figure \ref{fig1}.
\vspace{0.1cm}

\subsection{Numerical values} \label{numoriginal}
The following table sums up the numerical values considered in our setting:
\vspace{0.2cm}

\begin{center}
\begin{tabular}{|c|c|c|}
 \hline
\textbf{Variable} & \textbf{Value} & \textbf{Unit} \\
\hline
$d$ & $0.65$ & $m$\\
\hline
$C_{d}$ & $0.4$ & n.a\\
\hline
$T_{\max}$ & $5000$ & $N$\\
\hline
$C_{s}$ & $4.10^{-4}$ & $kg.s^{1}.N^{-1}$\\
\hline
$g$ & $9.81$ & $m.s^{-2}$\\
\hline
$\rho_{0}$ &  $1.225$ & $kg/m^{3}$\\
\hline
$h_{r}$ &  $7314$ & $m$\\
\hline
$u_{\max}$ & $2$ & n.a\\
\hline
$h_{c}$ & $250$ & $m$\\
\hline
$(k_{0}, k_{1})$ & (var, var) & $n.a$\\
 \hline
\end{tabular}
\end{center}
\vspace{0.2cm}

\noindent We consider the following boundary conditions:
\begin{subequations}  \label{boundcondbunt}
\begin{align}
&(x_{0}, h_{0}, v_{0}, \gamma_{0},m_{0})=(0 \text{m}, 0 \text{m}, 300 \text{m/s}, 80^{\circ}, 600 \text{kg}) \label{boundcondbunt1} \\
&(x_{f}, h_{f}, v_{f}, \gamma_{f}, m_{f})=(25000 \text{m}, 0 \text{m}, \text{*}, -80^{\circ}, \text{*}) \label{boundcondbunt2}
\end{align}
\end{subequations}
\begin{remark}[] \label{rem5}
As the thrust control $u_{1}(.)$ has been set to $1$, one has to let the final velocity $v_{f}$ free. If $u_{1}(.)$ was an active control, one could specify a desired value for $v_{f}$.
\end{remark}

\section{Optimal control in finite dimension} \label{Optcont}
The term finite dimension refers to the fact that the state vector $\xi(\cdot)$ belongs to the finite dimensional space $\mathbb{R}^{n}$. However, it is important to keep in mind that solving an optimal control problem like \eqref{contopt1} in finite dimension requires being able to solve an infinite dimensional optimization problem.
\subsection{General setting} \label{Genset}
In this section, we consider the following general control problem: given $\mathcal{M}_{0}$ and $\mathcal{M}_{1}$  two submanifolds of $\mathbb{R}^{n}$ we aim at controlling the autonomous nonlinear system:
\begin{equation}
\dot{\xi}(t)=f(\xi(t), \alpha(t)) \text{ on } [0, t_{f}]
\label{dyngen}
\end{equation}
while minimizing the cost:
\begin{equation}
C(t_{f}, \alpha):=\int_{0}^{t_{f}} f^{0}(\xi(t), \alpha(t)) dt+g(t_{f},\xi(t_{f}))
\label{costgen}
\end{equation}
and such that:
\begin{equation}
\xi(0) \in \mathcal{M}_{0} \text{ and } \xi(t_{f}) \in \mathcal{M}_{f}
\label{boundcond}
\end{equation}
\vspace{0.1cm}

\noindent In the above definition, we assume $f:\mathbb{R}^{n} \times \mathbb{R}^{p} \rightarrow \mathbb{R}^{n}$, $f^{0}: \mathbb{R}^{n} \times \mathbb{R}^{p} \rightarrow \mathbb{R}$
and $g:  \mathbb{R} \times \mathbb{R}^{n} \rightarrow \mathbb{R}$ are of class $C^{1}$. The control $\alpha(\cdot)$ belongs to the set $L^{\infty}([0, t_{f}], A)$ where $A$ is a subset of $\mathbb{R}^{p}$. In the following, $t_{f}$ can be free or fixed. We denote $\mathcal{A}$ the set of  \textit{admissible controls} i.e. such that the corresponding trajectories steer the system from an initial point of $\mathcal{M}_{0}$ to a final point of $\mathcal{M}_{f}$.
\vspace{0.1cm}

\begin{remark} \label{existence}
It is far from being obvious that there exists a solution to the previous optimal control problem. However, if one assumes, in addition to the above mentioned regularity assumptions on $f$, $f^{0}$ and $g$ that the following asumptions hold:
\begin{enumerate}
\item $A$ is a compact subset of $\mathbb{R}^{p}$
\item there exists $b>0$ such that any admissible trajectory $\xi_{\alpha}(\cdot)$ and the associated time $t_{f}$ are bounded by $b$:
\begin{equation}
\exists b>0 \hspace{0.1cm} \mid \hspace{0.1cm} \forall \alpha \in \mathcal{A}, \hspace{0.1cm} \forall t \in [0, t_{f}], \hspace{0.1cm} t_{f}+\|\xi_{\alpha}(t) \|_{\infty} \leqslant b
\label{boundtraj0}
\end{equation}
\item for any $\mu \in \mathbb{R}^{n}$, the set $V$ defined by:
\begin{equation}
V(\mu):=\left\{\left(f^{0}(\mu, a)+\gamma, f(\mu, a)\right)  \mid a \in A, \hspace{0.1cm} \gamma \geqslant 0 \right\}
\label{convV}
\end{equation}
is a convex subset of $\mathbb{R}^{n+1}$.
\end{enumerate}
\vspace{0.1cm}

\noindent Then there exists a solution $\alpha^{\star}(.)$ defined on an interval $[0, t(\alpha^{\star})]$ to the optimal control problem. The above mentioned assumptions ensure usual existence results, even if some of them could be weakened. For a survey of existence results in optimal control problems, refer to \cite{Ces1983} for instance.
\end{remark}
\vspace{0.1cm}

\begin{definition}[End point mapping]
Let $\alpha(.) \in\mathcal{A}$ be an admissible control. For given $t \in [0, t_{f}]$, and $\xi_{0} \in \mathbb{R}^{n}$, the end-point mapping $E_{t,\xi_{0}}$ is defined as below
\begin{subequations} \label{endpoint}
\begin{align}
E_{t, \xi_{0}}:&\mathcal{A} \longrightarrow \mathbb{R}^{n} \label{endpoint1} \\
&\alpha \mapsto \xi(t, \xi_{0}, \alpha) \label{endpoint2}
\end{align}
\end{subequations}
where $t \mapsto \xi(t, \xi_{0},\alpha)$ is the value at time $t$ of the solution to \eqref{dyngen} such that $\xi(0)=\xi_{0}$.
\end{definition}
\vspace{0.1cm}
\noindent It is well known that, if one endows $\mathcal{A}$ with the standard $L^{\infty}$ topology, then the end point mapping is $C^{1}$ on $\mathcal{A}$ and the optimal control problem \eqref{dyngen}-\eqref{costgen}-\eqref{boundcond} can be formulated in terms of end point mapping as follows:
\begin{equation}
\min \{C(t_{f}, \alpha), \hspace{0.1cm} \xi(0) \in \mathcal{M}_{0}, \hspace{0.1cm} E_{t_{f}, \xi_{0}}(\alpha) \in \mathcal{M}_{f}, \hspace{0.1cm} \alpha \in \mathcal{A} \}
\label{controptmapp}
\end{equation}
General considerations on optimal control in finite and infinite dimension can be found for instance in \cite{Tre2024}.
\vspace{0.2cm}

\noindent In the next paragraph, we give a statement of the Pontryagin maximum principle (PMP). It expresses necessary conditions for a pair $(\xi^{\star}(\cdot), \alpha^{\star}(\cdot))$ to be optimal for \eqref{dyngen}-\eqref{costgen}-\eqref{boundcond}.
\subsection{Pontryagin maximum principle} \label{PMP}
By the PMP (see \cite{Tre2024}, \cite{Pon1962}, \cite{Tre2008}), there exists $p^{0} \leq 0$ and an absolutely continuous mapping  $p^{\star}(.): [0, t_{f}] \rightarrow \mathbb{R}^{n}$, called \textit{adjoint vector},  such that, for almost every $t \in [0, t_{f}]$
\begin{subequations} \label{equapmp}
	\begin{align}
	&\dot{\xi}^{\star}(t) =\dfrac{\partial H}{\partial p}(\xi^{\star}(t),p^{\star}(t),p^{0},\alpha^{\star}(t)) \label{equapmp1} \\
	&\dot{p}^{\star}(t) =-\dfrac{\partial H}{\partial \xi}(\xi^{\star}(t),p^{\star}(t),p^{0},\alpha^{\star}(t))  \label{equapmp2}
	\end{align}
\end{subequations}
where the Hamiltonian $H$ is defined by:
\begin{equation}
H(\xi,p,p^{0},\alpha):=\langle p, f(\xi,\alpha) \rangle+p^{0}f^{0}(\xi,\alpha)
\label{Ham}
\end{equation} 
For almost every $t \in [0, t_{f}]$, the control $\alpha^{\star}$ maximizes the Hamiltonian $H$:
\begin{equation}
H(\xi^{\star}(t),p^{\star}(t),p^{0},\alpha^{\star}(t)) = \underset{v \in A}{\max}\hspace{0.1cm} H(\xi^{\star}(t),p^{\star}(t),p^{0},v)
\label{maxcond}
\end{equation} 
The adjoint vector satisfies the \textit{transversality conditions} (when the tangent space is well defined):
\begin{subequations} \label{transvp}
\begin{align}
&p^{\star}(0) \perp T_{\xi(0)} \mathcal{M}_{0} \label{transv1}\\
&p^{\star}(t_{f})-p^{0}\dfrac{\partial g}{\partial \xi}(t_{f}, \xi^{\star}(t_{f})) \perp  T_{\xi(t_{f})} \mathcal{M}_{f} \label{transv2}
\end{align}
\end{subequations} 
where $T_{\xi} \mathcal{M}$ is the notation for the tangent space of the submanifold $\mathcal{M}$ at $\xi$. If the final time $t_{f}$ is free, there is an additional transversality condition:
\begin{equation}
\underset{v \in A}{\max}\hspace{0.1cm} H(\xi^{\star}(t_{f}),p^{\star}((t_{f})),p^{0},v)=-p^{0}.\dfrac{\partial g}{\partial t}(t_{f}, \xi^{\star}(t_{f}))
\label{transvH}
\end{equation}
\begin{remark}[Autonomous case] \label{Hnull}
In the case where the final time is free one has:
\begin{equation}
\underset{v \in A}{\max}\hspace{0.1cm} H(\xi^{\star}(t),p^{\star}((t)),p^{0},v)=\mathrm{Cst}
\label{Hmaxnull}
\end{equation} 
\end{remark}
\vspace{0.1cm}

\noindent An \textit{extremal} of the optimal control problem is a fourth-tuple $(\xi(\cdot), p(\cdot), p^{0}, \alpha(\cdot))$ solution of \eqref{equapmp} and \eqref{maxcond}. If $p^{0}=0$, the extremal is said to be \textit{abnormal}, if $p^{0}<0$, then the extremal is said to be \textit{normal}.
\vspace{0.1cm}

\begin{remark}[Sufficient conditions for optimality] \label{Legendre}
We emphasize once more that the PMP gives a set of necessary conditions. For clarity, let us simplify the optimal control problem \eqref{dyngen}-\eqref{costgen}-\eqref{boundcond} by taking $A=\mathbb{R}^{p}$. If the \textit{strong Legendre condition} holds along a given extremal $(\xi(\cdot), p(\cdot), p^{0}, \alpha(\cdot))$, that is, there exists $r>0$ such that:
\begin{equation}
\dfrac{\partial^{2} H}{\partial \alpha^{2}}(\xi(\cdot), p(\cdot), p^{0}, \alpha(\cdot))(\nu, \nu) \leqslant -r. \|\nu\|^{2}, \hspace{0.2cm} \forall \nu \in \mathbb{R}^{p}
\label{legendrecond}
\end{equation} 
\end{remark}
\noindent then there exists $\epsilon> 0$ so that the trajectory $\xi(.)$ is locally optimal in $L^{\infty}$ topology on $[0, \epsilon]$. If the extremal is moreover normal, i.e. $p^{0} \neq 0$, then $\xi(\cdot)$ is locally optimal in $C^{0}$ topology on $[0, \epsilon]$. For more details on second order optimality conditions, refer for instance to \cite{AgrSac2004} or \cite{BonCaiTre2007}.
\vspace{0.1cm}

\begin{definition}[Singular arc] \label{singdef}
Assume $\mathcal{M}_{0}=\{\xi_{0}\}$. A control $\alpha_{s}(\cdot)$ defined on $[0, t_{f}]$ is said to be \textit{singular}
 if and only if the Fréchet differential $dE_{t_{f}, \xi_{0}}(\alpha_{s})$ is not surjective. The associated trajectory $\xi_{s}(\cdot)$ is called \textit{singular} trajectory.
\end{definition}
\vspace{0.1cm}

\noindent We recall the following two standard characterizations of singular controls (see \cite{BonChy2003}, \cite{Pon1962}). A control $\alpha_{s} \in  \mathcal{A}$ is singular if and only if the linearized system along the trajectory $\xi_{s}(\cdot, \xi_{0}, \alpha_{s})$ is not controllable. This is also equivalent to the existence of an absolutely continuous mapping $p_{s}: [0, t_{f}] \rightarrow \mathbb{R}^{n} \backslash \{0\}$ such that, for almost every $t \in[0,t_{f}]$,
\begin{subequations} \label{equasing}
	\begin{align}
	&\dot{\xi}_{s}(t) =\dfrac{\partial H_{s}}{\partial p}(\xi_{s}(t),p_{s}(t), \alpha_{s}(t)) \label{equasing1} \\
	&\dot{p}_{s}(t) =-\dfrac{\partial H_{s}}{\partial \xi}(\xi_{s}(t),p_{s}(t), \alpha_{s}(t))  \label{equasing2} \\
	&\dfrac{\partial H_{s}}{\partial \alpha}(\xi_{s}(t), p_{s}(t), \alpha_{s}(t))=0  \label{equasing3}
	\end{align}
\end{subequations}
where the Hamiltonian $H_{s}$ of the system is defined by:
\begin{equation}
H_{s}(\xi,p,\alpha):=\langle p, f(\xi,\alpha) \rangle
\label{Hamsing}
\end{equation} 
Note that singular trajectories coincide with projections of abnormal extremals for which the maximization condition \eqref{maxcond} reduces to $\partial H/\partial \alpha=0$. In case when the dynamics $f$ and the cost $f^{0}$ are linear in the control $\alpha$, a singular arc (restriction of an extremal to a subinterval $I$) corresponds to an arc along which one is unable to compute the control directly from the maximization condition of the PMP (at the contrary of bang-bang situation). Indeed, in this case, the above condition $\partial H/\partial \alpha=0$ along the arc means that some function (called switching function) vanishes identically along the arc. Then, it is well known that, in order to derive an expression of the control along such an arc, one has to differentiate this relation until the control appears explicitly. It is as well known that such singular arcs, whenever they occur, may be optimal. Their optimal status may be proved using generalized Legendre-Clebsch type conditions or the theory of conjugate points (see \cite{Rob1967}-\cite{Goh1966} or see \cite{AgrSac2004}-\cite{BonCaiTre2007} for a complete second-order optimality theory of singular arcs). 
\vspace{0.1cm}

\begin{remark} \label{singbound}
In case when $A$ is compact ($\alpha(\cdot)$ is constrained), the characterization \eqref{equasing} remains valid as soon as the control values remain in $\mathring{A}$.
\end{remark}

\section{Numerical methods in optimal control} \label{nummeth}
\subsection{Context}
Let us first recall that there are mainly two kinds of numerical approaches in optimal control: direct and indirect methods.
\vspace{0.1cm}

\noindent On the one hand, direct methods consist of discretizing the state and the control so as to reduce the problem to a constrained nonlinear optimization problem in finite dimension. The process is straightforward and it can be applied in systematic way to any optimal control problem. Another great advantage of the direct methods is that they do not require any prior information on the structure of the control. New variables or constraints can be added easily. However, the convergence may be difficult due to the large number of variables. 
\vspace{0.1cm}

\noindent On the other hand, indirect methods are based on the PMP which provides a set of necessary conditions for a local minimum. The problem is then reduced to a nonlinear system that is generally solved by a shooting method using a Newton-like algorithm. An advantage of such methods is their great accuracy when they converge. The drawbacks are the narrow radius of convergence and in some cases (presence of singular arcs, state constraints...) a prior theoretical work is necessary in order to identify the structure of the control. For more complete comparison between direct and indirect methods, refer to \cite{Tre2012}. The principles of both direct and indirect methods are recalled hereafter.
\subsection{Direct methods} \label{direct}
Direct methods consist of discretizing both the state and the control. After discretizing, the optimal control problem is reduced to a nonlinear optimization problem in finite dimension, or nonlinear programming problem, of the form:
\begin{equation} \label{NLP}
\underset{Z \in C}{\min} \hspace{0.1cm} F(Z)
\end{equation}
where $Z:=\left(\xi_{1},\xi_{2},...,\xi_{N},\alpha_{1},\alpha_{2},...,\alpha_{N}\right)$ and:
\begin{equation}
C:=\bigg\{Z \mid g_{i}(Z)=0, \hspace{0.1cm} i \in  [\![ 1,r]\!]; \hspace{0.3cm} g_{j}(Z) \leqslant 0, j \in [\![r+1,N]\!] \bigg\}
\label{constr}
\end{equation}
There exists an infinite number of variants, depending on the choice of finite-dimensional representations of the control and of the state, of the discretization of the extremal differential equations, and of the discretization of the cost functional. The discretization may be carried out in many ways, depending on the problem features. As an example, we may consider a subdivision $0=t_{0}<t_{1}<..<t_{N}=t_{f}$ of the interval $[0, t_{f}]$. We discretize the controls such that they are piecewise constant on this subdivision with values in $U$. Meanwhile, the differential equations may be discretized by an explicit Euler method, by setting $h_{i}:=t_{i+1}-t_{i}$, we get $\xi_{i+1}=\xi_{i}+h_{i}.f(\xi_{i}, \alpha_{i})$. The cost may be discretized by a quadrature procedure.  We refer to \cite{Betts2010} for a thorough description of many direct approaches in optimal control. 
\vspace{0.1cm}

\noindent Then, to solve the optimization problem \eqref{NLP} under the constraints \eqref{constr}, there is also a large number of possible methods: gradient methods, penalization, quad-Newton, dual methods... We refer the reader to any good textbook of numerical optimization. In the frame of our study, we use the optimization routine IpOpt combined with the automatic differentiation code AMPL on a standard desktop machine. For more information, refer to \cite{FouGayKer2002}-\cite{WacBie2006}. Alternative variants of direct methods are the collocation methods, the spectral or pseudo-spectral methods, the probabilistic approaches, etc.
\vspace{0.2cm}

\noindent Another approach to optimal control problems that can be considered as a direct method, consists in solving the Hamilton-Jacobi-Bellman equation satisfied (in the viscosity sense) by the value function which is of the form:
\begin{equation}
\dfrac{\partial v}{\partial t}+H\left(\xi, \dfrac{\partial v}{\partial \xi}\right)=0
\end{equation}
The value function $v$ is the optimal cost for the optimal control problem starting from a given point $(t, \xi)$. Once the value function estimated, optimal controls can be deduced from the Dynamic Programming Principle.
\subsection{Indirect methods} \label{indirect}
In indirect approaches, the PMP (first order necessary conditions for optimality) is applied to the optimal control problem. This reduces the problem to a nonlinear system of $n$ equations with $n$ unknowns generally solved by Newton-like methods. The indirect methods are also called shooting methods. The principle of the simple shooting and of the multiple shooting method is recalled herafter.
\vspace{0.2cm}

\noindent \textbf{Simple shooting method}: By setting $z(t):=(\xi(t), p(t))$, assume that, by using \eqref{maxcond}, the optimal control can be expressed as a function of the state and the adjoint variable. Then the extremal system \eqref{equapmp1}-\eqref{equapmp2} can be written under the form: $\dot{z}(t)=F(z(t))$. The initial and final conditions \eqref{boundcond}, the transversality conditions on the adjoint \eqref{transvp} and Hamiltonian \eqref{transvH} can be written as $R(z(0),z(t_{f}),t_{f})=0$. Finally, we obtain a two-point boundary value problem
\begin{equation}
\dot{z}(t)=F(t, z(t)), \hspace{0.2cm} R(z(0), z(t_{f}), t_{f})=0
\label{FR}
\end{equation}
Let $z(t, z_{0})$ be the solution of the Cauchy problem:
\begin{equation}
\dot{z}(t)=F(t, z(t)), \hspace{0.2cm} z(0)=z_{0}
\label{CP}
\end{equation}
The two-point boundary problem \eqref{FR} consists in finding a zero of the equation:
\begin{equation}
R(z(0), z(t_{f}, z_{0}), t_{f})=0
\label{eqR}
\end{equation}
This problem can be solved by Newton-like or other other iterative methods.
\vspace{0.2cm}

\noindent \textbf{Multiple shooting method}: The drawback of the single shooting method is the sensitivity of the Cauchy problem to the initial condition $z_{0}$. The multiple shooting aims at a better numerical stability by dividing the time interval $[0, t_{f}]$ into $N$ subintervals $[t_{i}, t_{i+1}]$ and considering as unknowns the values of $z_{i}=(\xi(t_{i}), p(t_{i}))$ at the beginning of each subinterval. The application of the PMP to the optimal control problem yields a multi-point boundary value problem, which consists in finding $Z=(p(0), t_{f}, z_{i})$ for $i=1.., N-1$ such that the differential equation:
\begin{equation}
\dot{z}_{i}(t)=F(t, z(t))=\left\{
\begin{array}{l}
F_{0}(t, z(t)), \hspace{0.5cm} t_{0} \leqslant t \leqslant t_{1} \\
F_{1}(t, z(t)), \hspace{0.5cm} t_{1} \leqslant t \leqslant t_{2} \\
...., \\
F_{N-1}(t, z(t)), \hspace{0.1cm} t_{N-1} \leqslant t \leqslant t_{N}
\end{array}
\right.
\end{equation}
and the constraints:
\begin{subequations}
\begin{align*}
&\xi(0) \in \mathcal{M}_{0}, \hspace{0.1cm} \xi(t_{f}) \in \mathcal{M}_{f},  \hspace{0.1cm}  p^{\star}(0) \perp T_{\xi(0)} \mathcal{M}_{0} \\
&p^{\star}(t_{f})-p^{0}\dfrac{\partial g}{\partial \xi}(t_{f}, \xi^{\star}(t_{f})) \perp  T_{\xi(t_{f})} \mathcal{M}_{f} \\
& z_{i}(t_{i}^{-})=z_{i}(t_{i}^{+}), \hspace{0.1cm} i=1,...,N-1, \hspace{0.1cm} H(t_{f})=0
\end{align*}
\end{subequations}
are satisfied. The nodes of the multiple shooting may involve the switching times (at which the switching function changes sign) and the junction times (entry, contact or exit times) with boundary arcs for instance. In this case, an a priori knowledge of the solution structure is required. For more details on the shooting methods, see \cite{StoBul2002}.

\subsection{Methods implemented in the paper} \label{nummethodpaper}
Even though the direct methods will be implemented and evaluated, the shooting method based on the application of the PMP will be privileged in this paper. The reasons for this choice are the following: the first one is the enhanced accuracy which is crucial in the applications of the optimal control to the aerospace field. Another reason is the possibility to adapt this method to a specific problem: indeed, if adequately initialized, the shooting method can provide a quasi instantaneous solution compliant with the real time implementation. Moreover, in case of deviation with respect to the optimal solution, one is able to compute online a new one. These requirements are important for the applications in the aerospace field and cannot, in general, be satisfied by direct methods as soon as the problem becomes too much complex. However, as explained previously the initialization is in general a hard and intricate task, which requires a detailed study of the optimal problem structure and sometimes additional mathematical tools coming for instance from the geometric optimal control (refer for instance to \cite{Tre2012}). 
\vspace{0.2cm}

\noindent As it is well known, the main issue of the shooting methods based on the Newton-like method is the initialization: indeed, due to the small convergence radius, the first guess shall be "close enough" to the solution. In recent years, the numerical continuation has become a powerful tool to overcome this difficulty.  Section \ref{contmeth} recalls some mathematical concepts of the continuation approaches, with a focus on the numerical implementation of these methods.

\section{Continuation methods} \label{contmeth}
\subsection{Existence results and discrete continuation} \label{existcont}
The basic idea of continuation (also called homotopy) method is to solve a difficult problem step by step starting from a simpler problem by parameter deformation. The theory and practice of the continuation methods are well-spread (see, e.g., \cite{AllGeo1990}, \cite{Rhei2000}, \cite{Wat2001}). Combined with the shooting method derived from the PMP, a continuation methods consists in deforming the problem into a simpler one (that can be easily solved) and then solving a series of shooting problems step by step to ``come back'' to the original problem. 
\vspace{0.1cm}

\noindent One difficulty of the homotopy methods lies in the choice of a sufficiently regular deformation that allows the convergence along the path. The starting problem should be easy to solve and the path between this starting problem and the original problem should be easy to model. Another difficulty is to numerically follow the path between the starting problem the original problem. This path is parametrized by a parameter denoted usually by $\lambda$.
\vspace{0.1cm}

\noindent The choice of the homotopic parameter may require considerable physical insight into the problem. This parameter may be defined either artificially according to some intuition, or naturally by choosing physical parameters of the system, or by combination of both.
\vspace{0.1cm}
	
\noindent Assume one has to solve a system of $N$ nonlinear equations in $N$ dimensional variable $Z$:
\begin{equation}
F(Z)=0
\label{homoF}
\end{equation}
where $F: \mathbb{R}^{n} \longrightarrow \mathbb{R}^{n}$ is a smooth map. We define a deformation $G: \Omega \times [0, 1] \longrightarrow \mathbb{R}^{n}$ such that:
\begin{equation}
G(Z,0)=G_{0}(Z), \hspace{0.2cm} G(Z,1)=F(Z)
\label{homodef}
\end{equation}
where $G_{0}: \mathbb{R}^{n} \longrightarrow \mathbb{R}^{n}$ is a smooth map having known zeros.
\vspace{0.1cm}

\noindent A zero path is a curve $c(s) \in G^{-1}(0)$ where $s$ represents the arc length. We would like to trace a zero path starting from a point $Z_{0}$ such that $G(Z_{0},0)=0$ and ending at a point $Z_{f}$ such that $G(Z_{f},1)=0$ (see figure \ref{fig4}).

\begin{figure}[h]
\begin{center}
\captionsetup{justification=centering}
\includegraphics[width=0.5\linewidth]{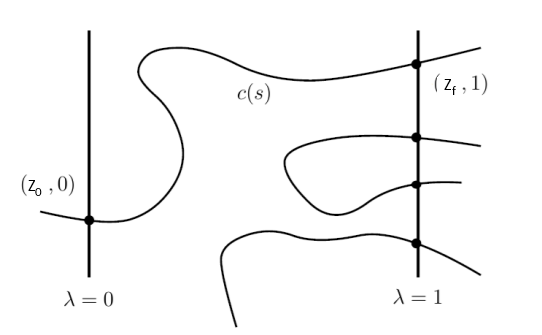}
\caption{Different zero paths}
\label{fig4}	
\end{center}
\end{figure}
\FloatBarrier

\noindent The first question to address is the existence of zero paths, since the feasability of the continuation method relies on this assumption. Actually, it has been proved in \cite{Tre2012}-\cite{TreZhuCer2017}, that the local feasability of the continuation method is closely related to the three following conditions:
\begin{enumerate}
\item there are no minimizing abnormal extremals
\item there are no minimizing singular controls, meaning the mapping $dE_{t_{f}, \xi_{0}}(\alpha)$ is surjective (see definition \ref{singdef}).  
\item there are no conjugate points: for precise definition of conjugate points, refer to \cite{Tre2012}-\cite{TreZhuCer2017}. Actually the absence of conjugate points can be numerically tested (see \cite{BonCaiTre2007}).
\end{enumerate}
Finally, it has to be noticed that, despite of local feasability, the zero paths may not be globally defined for any $\lambda \in [0, 1]$. The path could cross some singularity or diverge to infinity before reaching $\lambda=1$. The first possibility can be discarded by assuming (2) and (3) over all the domain $\Omega$ and for every $\lambda \in [0, 1].$ The second possibility is related to some properness properties of the exponential mapping (see \cite{BonTre2001}-\cite{Tre2000}).
\subsection{Numerical tracking of the zero paths}
There exists many numerical algorithms to track a zero path. Among these algorithms, the simplest one is called the \textit{simple continuation} procedure.
\vspace{0.2cm}

\noindent \textbf{Basic continuation}: The continuation parameter $\lambda$ is discretized by $0=\lambda_{0}<\lambda_{1}<..<\lambda_{n}=1$ and the sequence of problems $G(Z,\lambda_{i})=0$ is solved to end up with a zero point of $F(Z)$. If the increment $\Delta \lambda=\lambda_{i+1}-\lambda_{i}$ is small enough, then the solution $Z_{i}$ associated to $\lambda_{i}$ to $G(Z_{i},\lambda_{i})=0$ is generally close to the solution of $G(Z,\lambda_{i+1})=0$. An implicit assumption is made that we are able to compute the solution at the first step of the continuation procedure, namely a zero $Z_{0}$ for $G(Z, 0)=0$. The simple continuation algorithm is detailed below:
\begin{algorithm}[H]
\caption{Simple continuation procedure}
\begin{algorithmic}
\State  \textbf{Initialization:} $\lambda=0$, $Z=Z_{0}$, $\Delta \lambda \in (\Delta \lambda_{\min}, \Delta \lambda_{\max})$
\While{$\Delta \lambda \in (\Delta \lambda_{\min}, \Delta \lambda_{\max})$ and $\lambda \leqslant 1$}
\State $\Delta \lambda=\min(\Delta \lambda, 1- \lambda)$
\State $\tilde{\lambda}=\lambda+\Delta \lambda$
\State Look for $\widetilde{Z}$  zero of $G(\widetilde{Z}, \tilde{\lambda})=0$
\If{success}
\State $Z =\widetilde{Z}$
\State $\lambda =\tilde{\lambda}$
\State $\Delta \lambda=2. \Delta \lambda$
\Else
\State $\lambda=\tilde{\lambda}-\Delta \lambda$
\State $\Delta \lambda =\Delta \lambda/2$
\EndIf
\EndWhile
\If{success}
\State The continuation procedure is successful
\Else 
\State The continuation procedure has failed
\EndIf
\end{algorithmic}
\end{algorithm}
\noindent \textbf{Continuation with linear prediction}: Behind this procedure is the idea that we can do better than just using $Z_{\lambda}$ to approximate $Z_{\lambda+\Delta \lambda}$. Assume that we have already made two resolutions, yielding $Z_{\lambda_{1}}$ and $Z_{\lambda}$,  
for two values $\lambda_{1}$ and $\lambda$ such that $\lambda_{1}< \lambda$. Assuming some regularity on the path of zeros, an approximation
of $Z_{\lambda+ \Delta \lambda}$ for a new value $\lambda+\Delta \lambda$ is given by
\begin{equation}
Z_{\lambda+\Delta \lambda}=Z_{\lambda}+\dfrac{\Delta \lambda}{\lambda-\lambda_{1}}\left(Z_{\lambda}-Z_{\lambda_{1}} \right)
\label{contpredlin}
\end{equation}
\begin{figure}
\centerline{\includegraphics[width=0.5\linewidth]{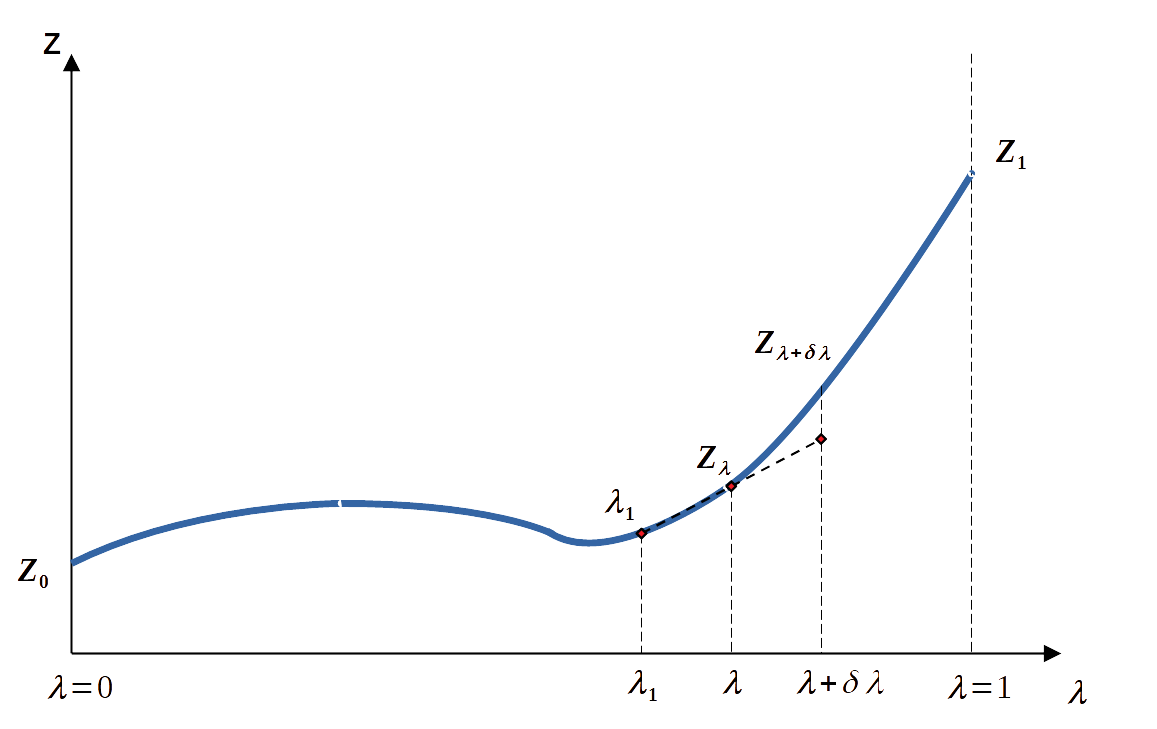}}
\caption{Continuation with linear prediction}
\label{fig5}	
\end{figure}
\FloatBarrier
\noindent This is the procedure implemented throughout this article each time a continuation is performed, and one experimentally witnesses an improvement in the runtime of the algorithm.

\section{The analysis of the $(\textbf{OCP})_{k_{0},k_{1}}$ control problem}
In the equations given below, the time dependence has been skipped on purpose for the sake of readability. We use herafter the notations of the previous section. Here, the control is in dimension 1 and the set $U=[-u_{\max}, u_{\max}]$.

\subsection{Computation of the extremals}
We apply the Pontryagin maximum principle recalled in Section \ref{PMP}. By denoting $\langle , \rangle$ the usual scalar product in $\mathbb{R}^{5}$ and the adjoint vector by:
\begin{equation}
p=\left(\begin{array}{c}
p_{x} \\
p_{h} \\
p_{v} \\
p_{\gamma} \\
p_{m}
\end{array}
\right)
\end{equation}
the Hamiltonian of the optimal control problem \eqref{contopt1} is given by:
\begin{align*}
H_{0}&:=\langle p, f(\xi, u) \rangle + p^{0}f^{0}(\xi, u) \\
&= p_{x}v\cos \gamma+p_{h}v\sin \gamma+p_{v}\bigg(\dfrac{T_{\max}\left(1+C_{s}v\right)-D(h,v)}{m}-g\sin \gamma\bigg) \\
&\hspace{2cm}+\dfrac{p_{\gamma}}{v}\left(\dfrac{L(h,v,u)}{m}-g\cos \gamma \right)-C_{s}T_{\max}p_{m}+p^{0}\bigg(k_{0}+k_{1}\left(\dfrac{h-h_{c}}{h_{c}} \right)^{2}\bigg)
\end{align*}
We assume hereafter that the optimal state is not abnormal thus we take $p^{0}=-1$. The adjoint equations can then be derived and one obtains:
\begin{subequations} \label{adjeqsing}
\begin{align} 
&\dot{p}_{x}=0 \label{adjeqsing1} \\
&\dot{p}_{h}=\dfrac{p_{v}}{m}\dfrac{\partial D}{\partial h}(h,v)-\dfrac{p_{\gamma}}{mv}\dfrac{\partial L}{\partial h}(h,v,u)+\dfrac{2k_{1}}{h_{c}^{2}}(h-h_{c}) \label{adjeqsing2} \\
&\dot{p}_{v}=-p_{x}\cos \gamma-p_{h}\sin \gamma-\dfrac{p_{v}}{m}C_{s}T_{\max}+\dfrac{p_{v}}{m}\dfrac{\partial D}{\partial v}(h,v)-\dfrac{p_{\gamma}}{m}\dfrac{\partial(L/v)}{\partial v}-\dfrac{p_{\gamma}}{v^{2}}g\cos \gamma \label{adjeqsing3} \\
&\dot{p}_{\gamma}=p_{x}v\sin \gamma-p_{h}v\cos \gamma+p_{v}g\cos \gamma-\dfrac{p_{\gamma}}{v}g\sin \gamma \label{adjeqsing4}\\
&\dot{p}_{m}=\dfrac{p_{v}}{m^{2}}\bigg(T_{\max}(1+C_{s}v)-D(h,v)\bigg)+\dfrac{p_{\gamma}}{vm^{2}}L(h,v,u)  \label{adjeqsing5}
\end{align}
\end{subequations}
The maximization condition \eqref{maxcond} reads:
\begin{equation}
u \in \underset{|r| \leqslant u_{\max}}{\arg \max }\left(p_{\gamma}r\right)
\end{equation}
which leads to:
\begin{equation}
u(t)=\left\{
\begin{array}{l}
u_{\max}\dfrac{p_{\gamma}(t)}{|p_{\gamma}(t)|} \hspace{1.1cm} \text{ if } p_{\gamma}(t) \neq 0 \\
\text{to be determined } \hspace{0.3cm} \text{if}\hspace{0.1cm} p_{\gamma}(t)=0
\end{array}
\right.
\label{optcontlin}
\end{equation}
Moreover, the adjoint transversality conditions \eqref{transvp} read:
\begin{equation}
p_{v}(t_{f})=0 \text{ and } p_{m}(t_{f})=0
\end{equation}
Finally in our case (no terminal cost, free final time and autonomous dynamics and cost), \eqref{transvH} leads to:
\begin{equation}
H_{0}(\xi(t_{f}),p(t_{f}),-1,u(t_{f}))=0
\label{Ham0tf}
\end{equation}
\subsection{Overview of the optimal control structure} \label{sing_control}
The control problem $(OCP)_{k_{0},k_{1}}$ consists in steering the single input control-affine system from $\xi_{0} \in \mathbb{R}^{5}$ to the final target $\xi_{f} \in \mathbb{R}^{3}$ minimizing the cost \eqref{contopt1}. 
In \eqref{optcontlin} the adjoint variable $p_{\gamma}$ is classically called the \textit{switching function}. As explained in the previous section, the structure of $u(.)$ may include singular arcs, which occur whenever the adjoint state $p_{\gamma}$ vanishes along a subinterval of $[0, t_{f}]$. Of course, this does not necessarily happen. Note that when the switching function does not vanish over a subinterval, then the arc is either "bang max" (i.e. $u(t)=u_{\max}$) or "bang min" (i.e. $u(t)=-u_{\max}$). Briefly, the optimal policy may be a  concatenation of bang and/or singular arcs.
\newline
If $I$ is a subinterval of $[0, t_{f}]$ where a singular arc occurs, one has:
\begin{equation}
p_{\gamma}(t)=0 \text{ for any } t \in I
\label{pgzero}
\end{equation}
The usual method to compute the associated singular control is to differentiate repeatedly the switching function until the control appears explicitely. More precisely, from \eqref{pgzero} we have as well:
\begin{equation}
\dot{p}_{\gamma}(t)=0 \text{  a.e. on } I
\label{dpgdtzero}
\end{equation}
and differentiating once more, we obtain
\begin{equation}
\ddot{p}_{\gamma}(t)=0 \text{  a.e. on } I
\label{d2pgdt2zero}
\end{equation}
The singular control $u_{s}$ is expected to appear explicitly in the latter relation. The computations lead to an equation of the form:
\begin{equation}
A(\xi(t),p(t)).u_{s}(t)=B(\xi(t),p(t))
\label{alfasingeq}
\end{equation}
almost everywhere on $I$. $u_{s}$ can be deduced from the last equation as this relation should be generically nontrivial, that is, the coefficient $A$ should not be equal to zero. For more details, refer to \cite{BonMarTre2008} for instance.
\subsection{Singular control in the guidance problem} \label{sing_control_guidance}
In our case, $u_{s}$ can be analytically calculated:
\begin{equation}
u_{s}(t)=\dfrac{mv}{\bar{q}(h,v)S}\left(\dfrac{g\cos\gamma}{v}+\dfrac{A(\xi,p)}{B(\xi,p)}\right)
\label{sing_contr_guidance}
\end{equation}
where:
\begin{align*}
A(\xi,p)&=(p_{h}\cos \gamma-p_{x}\sin \gamma)\left(\dfrac{T_{\max}(1+C_{s}v)-D(h,v)}{m}-g\sin \gamma \right) \\
&\hspace{1.1cm}+g\cos \gamma\left(p_{x}\cos \gamma+p_{h}\sin \gamma+\dfrac{p_{v}C_{s}T_{\max}}{m}-\dfrac{p_{v}}{m} \dfrac{\partial D}{\partial v}\right) \\
&\hspace{3.1cm}+v\cos \gamma \left(-\dfrac{p_{v}D(h,v)}{mh_{r}} +\dfrac{2k_{1}}{h_{c}^2}(h-h_{c})\right)
\end{align*}
and
\begin{align*}
B(\xi,p)&=v(p_{h}\sin \gamma+p_{x}\cos \gamma)-gp_{v}\sin \gamma
\end{align*}
We recall that the expression \eqref{sing_contr_guidance} is valid as long as the switching function vanishes i.e $p_{\gamma}(t)=0$. \\
\newline
In practice, the main difficulty is to infer the structure of the optimal control meaning the chronological sequence of bang/singular (if any) arcs. An efficient way to do so is to run a direct method.
\subsection{The direct method approach} \label{ampl_ex}
We consider the case where $k_{0}=k_{1}=1$. We compute optimal solution of \eqref{contopt1} for different values of cruise altitude $h_{c}$ and boundary conditions. The table below summarizes the values that has been changed with respect to initial boundary conditions \eqref{boundcondbunt}:
\vspace{0.2cm}
\begin{center}
\begin{tabular}{|c|c|c|c|}
 \hline
\textbf{Case} & $\bm{\gamma_{0}}$ & $\bm{\gamma_{f}}$ & $\bm{h_{c}}$ \\
\hline
$1$ & $45^{\circ}$ & $-45^{\circ}$ & $1500m$\\
\hline
$2$ & $80^{\circ}$ & $-80^{\circ}$ & $250m$\\
\hline
$3$ & $15^{\circ}$ & $-15^{\circ}$ & $500m$\\
\hline
\end{tabular}
\end{center}
\vspace{0.2cm}
To implement the direct method, we use the optimization routine IpOpt (see  \cite{WacBie2006}) combined with the automatic differentiation code AMPL (see \cite{FouGayKer2002}) on a standard desktop machine. We discretize \eqref{contopt1} using a simple implicit Crank Nicholson method with $N=1000$ time steps. Some outputs are presented  hereafter:
\begin{figure}[h]
\begin{center}
\captionsetup{justification=centering}
\includegraphics[width=1.0\linewidth]{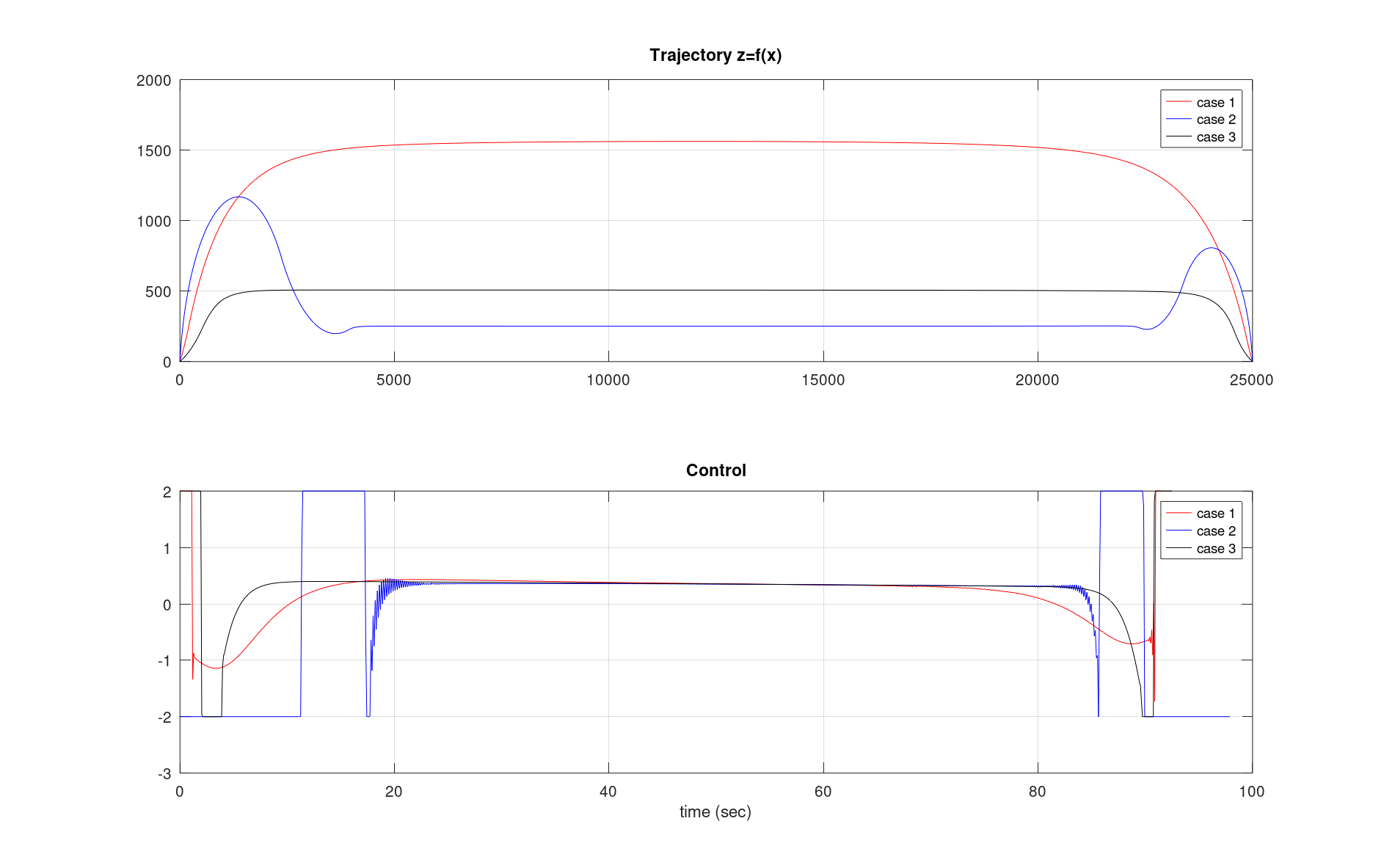}
\caption{Optimal control structure: 3 cases}
\label{fig5}	
\end{center}
\end{figure}
\FloatBarrier
\subsection{Formal analysis}
At the light of the observed results, one can infer the following:
\begin{itemize}
\item the cruise altitude penalization in the cost seems to be efficient in order to ensure the level flight at the specified altitude with acceptable error. This observation validates the numerical choices for $k_{0}$ and $k_{1}$ parameters.
\item during the level flight, one has $h(t) \approx h_{c}$ and $\gamma(t) \approx 0$. Actually, along the central arc, some coordinates of the optimal state, namely $h$ and $\gamma$ remain quasi constant: one is in presence of a partial turnpike phenomenon. The turnpike phenomenon has been widely described and analyzed in the literature (see Section \ref{sec2} for introduction and references). The partial turnpike is a variant of the latter which involves some coordinates of the state but not all of them (see \cite{Tre2023} for an example). This observation is crucial with respect to the continuation strategy we implement in section \ref{contphase2}.
\item the level flight arc at cruise altitude is likely singular (however we do not prove this assumption).
\item depending on the boundary conditions, there may be transitory initial (resp. final) ``bang'' arcs.
\end{itemize}

\noindent In the sequel of the paper we fix $k_{0}=k_{1}=1$.
\vspace{0.2cm}

\noindent In order to implement the shooting method based on the PMP, one is facing with two challenges:
\begin{enumerate}
\item the optimal structure of the control is not clearly identified particularly near $t=0$ and $t=t_{f}$ (as it depends strongly on the boundary conditions).
\item due to the narrow radius of convergence, the dimension of the model ($5$) makes the initialization intricate.
\end{enumerate}

\noindent To overcome the first difficulty, a classical method consists in regularizing the original problem by adding a quadratic term either to the dynamics or to the cost. This method has been widely used in the literature: in \cite{BonMarTre2008}, in the Goddard’s problem, the authors implement the homotopic approach based on the quadratic regularization of the cost in order to tackle with the nonsmoothness of the optimal control. In \cite{HabMarGer2004}, the maximum mass orbital transfer problem for a low thrust propulsion system is solved. In particular, the authors discover the complex structure of the optimal control (with multiple bang arcs) without any a priori assumptions thanks to the identical strategy. Note that this type of approach has been as well adopted in \cite{Mar2005} or \cite{MarGer2007}.
\vspace{0.2cm}

\noindent Consequently, in what follows we consider the regularized cost:
\begin{equation}
J_{k}(t_{f},u):=\displaystyle{\int_{0}^{t_{f}}}f_{k}^{0}(\xi(s), u(s)) ds:=\displaystyle{\int_{0}^{t_{f}}}\left(1+\dfrac{(h(s)-h_{c})^2}{h_{c}^2}+k.u(s)^{2}\right)ds
\label{costreg}
\end{equation}
where $k \geqslant 0$ is a regularization parameter, penalizing the $L^{2}$-norm of the control $u$. 
\vspace{0.2cm}

\noindent From now on, considering the optimal control problem \eqref{contopt1} we denote:
\begin{itemize}
\item $(\textbf{OCP})$ the case where $k_{0}=k_{1}=1$.
\item $(\textbf{OCP})_{k}$ the case where one considers the regularized cost \eqref{costreg}.
\end{itemize}

\noindent To overcome the second difficulty, we are going the define a simpler problem (called Dubins Fuller Regularized Problem $(\textbf{DFRP})_{k}$) in dimension $3$ for which one is able to initialize easily the shooting method. Then we are going to deform it into $(\textbf{OCP})_{k}$ by performing continuations over the dynamics and the boundary conditions. 
\subsection{Strategy for the shooting method}
In order to solve numerically $(\textbf{OCP})$ one proceeds as follows:
\begin{enumerate} 
\item Numerical solving of $(\textbf{DFRP})_{k}$ and continuations on the dynamics and boundary conditions leading to $(\textbf{OCP})_{k}$
\item Decreasing continuation on the regularization parameter $k$ of $(\textbf{OCP})_{k}$ in order to tend towards the initial problem $(\textbf{OCP})$ 
\end{enumerate} 
The following scheme illustrates the principles mentioned above:
\begin{figure}[h]
\begin{center}
\captionsetup{justification=centering}
\includegraphics[width=0.8\linewidth]{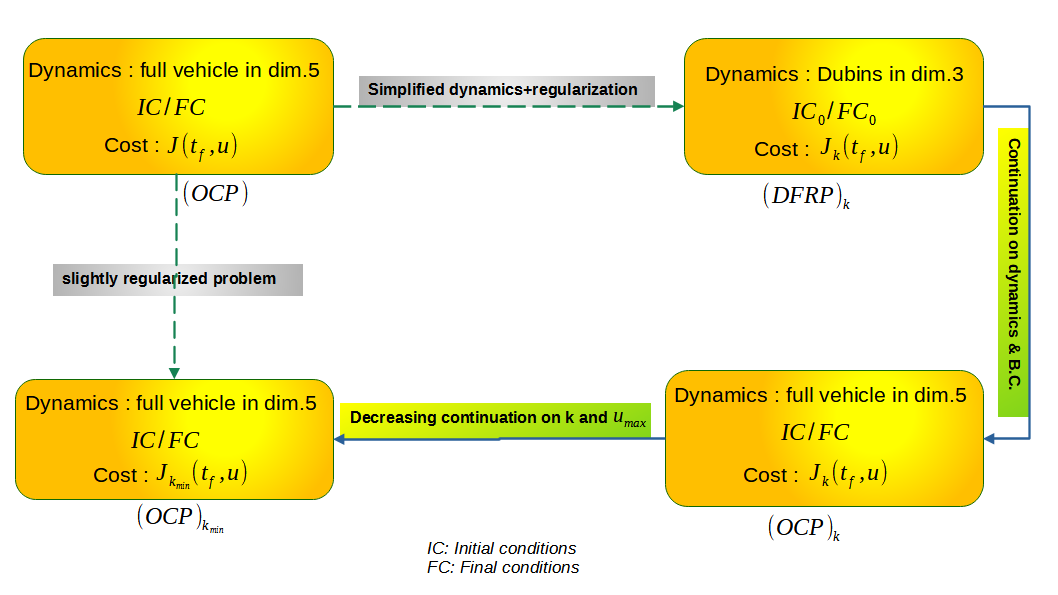}
\caption{Principles of numerical solving}
\label{fig6}	
\end{center}
\end{figure}
\FloatBarrier
\noindent \hspace{0.1cm}

\section{From $(\textbf{DFRP})_{k}$ to $(\textbf{OCP})_{k}$} \label{contphase1}
\subsection{Definition of $(\textbf{DFRP})_{k}$}
Let us consider the following kinematic equations in the vertical plane:
\begin{subequations} \label{dubdyn}
\begin{align}
&\dot{x} =v \cos \gamma \label{dubdyn1} \\
&\dot{h}=  v \sin \gamma \label{dubdyn2} \\
&\dot{\gamma}= u \label{dubdyn3}
\end{align}
\end{subequations}
where $v>0$ is the longitudinal speed, assumed to be fixed. $(x,h, \gamma)$ are the state variables defining respectively the position and the pitch angle of the vehicle. $u$ designates the pitch angular rate, which enables to maneuver and thus will be considered as a control variable. The rate is  limited and $|u| \leqslant u_{\max}$.
\vspace{0.2cm}

\noindent Actually, the \emph{Dubins} dynamics \eqref{dubdyn} is an elementary version of the vehicle model  as it can be obtained from the dynamics given in \eqref{contopt1} by making the following assumptions:
\begin{itemize}
\item neglect the gravity $g$ and assume the air density $\rho(\cdot)$ constant function of height $h$
\item consider the speed $v$ and the mass $m$ remain constant (at first order) during the flight. This can be achieved by setting $\dot{v}=0$ and $C_{s}=0$ in \eqref{contopt1}. Notice that these assumptions are far from being unrealistic in the case where the covered speed range is relatively small and the mass variation limited provided that the time interval remains reasonable.
\end{itemize}
By denoting the reduced dynamics $f_{red}(\mu,u):=\left(\begin{array}{ccc}
		v\cos \gamma \\
		v\sin \gamma \\
		u
	\end{array}
	\right)$, \eqref{dubdyn} can be written as $\dot{\mu}=f_{red}(\mu, u)$, with $\mu(0)=\mu_{0}$ and $\mu(t_{f})=\mu_{f}$.
\vspace{0.2cm}

\noindent $(\textbf{DFRP})_{k}$ consists in minimizing the cost $J_{k}(t_{f},u)$  under the dynamic constraint \eqref{dubdyn} while going from an initial state $\mu_{0}=(x_{0},h_{0}, \gamma_{0})$ to a final state $\mu_{f}=(x_{f},h_{f}, \gamma_{f})$.
\vspace{0.2cm}

\begin{remark} \label{remreg} The regularization of this simpler control problem is motivated by the fact that, as in the initial optimal control problem, the Hamiltonian depends linearly on the control.
\end{remark}
\begin{remark} \label{remdfrp} The acronym  \textbf{DFRP} designates the ``\textbf{D}ubins \textbf{F}uller \textbf{R}egularized \textbf{P}roblem'' which is the regularized version of the so called ``\textbf{D}ubins \textbf{F}uller'' problem (where $k=0$). For more details see \cite{DF2024}, where the model has been introduced and analyzed.
\end{remark}
\subsection{Parametrized dynamics}
Let us consider the following 5D dynamical model: $\dot{\xi}=f_{\lambda}(\xi, u)$ where $f_{\lambda}: \mathbb{R}^{5} \longrightarrow  \mathbb{R}^{5}$ is given by:
\begin{equation}
f_{\lambda}(\xi,u):=\left(\begin{array}{c}
v\cos \gamma \\
v\sin \gamma  \\
\lambda\left(\dfrac{T_{\max}(1+C_{s}v)-D(h,v)}{m} - g\sin \gamma \right)  \\
a(\xi,u).\lambda+u \\
-\lambda C_{s}T_{\max}
\end{array} \right)
\label{flambda}
\end{equation}
with $a(\xi,u):=\left(\dfrac{\bar{q}(h,v)S}{mv}-1\right)u-g\dfrac{\cos \gamma}{v}$. The dynamics is parametrized by the parameter $\lambda \in [0, 1]$.
\vspace{0.1cm}

\noindent For $\lambda=0$, one has the simplified dynamics \eqref{dubdyn} whereas for $\lambda=1.0$, one has the full dynamics as defined in \eqref{contopt1}. Consequently, it seems natural to perform an increasing continuation over $\lambda$ in order to ``connect'' both dynamics.
\subsection{Equations arising from the PMP}
The Hamiltonian of the regularized problem associated to the parametrized dynamics \eqref{flambda} is given by:
\begin{equation} 
H_{\lambda}^{k}(\xi,p, p^{0}, u):=\langle p, f_{\lambda}(\xi, u) \rangle + p^{0}f_{k}^{0}(\xi, u)
\label{Ham5dreg}
\end{equation}
We maintain the assumption of normality here by taking $p^{0}=-1$. 
\vspace{0.1cm}

\noindent The adjoint equations $\dot{p}=-\dfrac{\partial H_{\lambda}^{k}}{\partial \xi}$ read along each coordinate:
\begin{subequations} \label{adjeqparamdyn}
\begin{align*} 
&\dot{p}_{x}=0 \\
&\dot{p}_{h}=\lambda \left(\dfrac{p_{v}}{m} \dfrac{\partial D}{\partial h}(h,v)-\dfrac{p_{\gamma}}{m v}\dfrac{\partial L}{\partial h}(h,v,u)\right)+\dfrac{2}{h_{c}^{2}}(h-h_{c}) \\
&\dot{p}_{v}=-p_{x}\cos \gamma-p_{h}\sin \gamma+\lambda \left(\dfrac{p_{v}}{m}\dfrac{\partial D}{\partial v}(h,v)-\dfrac{p_{\gamma}}{m} \dfrac{\partial(L/v)}{\partial v}-\dfrac{p_{\gamma}}{v^{2}} g \cos \gamma-\dfrac{p_{v}}{m} C_{s} T_{\max} \right)  \\
&\dot{p}_{\gamma}=p_{x} v \sin \gamma-p_{h} v \cos \gamma+\lambda p_{v} g \cos \gamma-\lambda \dfrac{p_{\gamma}}{v} g \sin \gamma  \\
&\dot{p}_{m}=\lambda \dfrac{p_{v}}{m^{2}}\bigg(T_{\max}(1+C_{s}v)-D(h,v)\bigg)+\lambda\dfrac{p_{\gamma}}{v m^{2}}L(h,v,u)  
\end{align*}
\end{subequations}
Following the maximization condition \eqref{maxcond}, the optimal control $u$ is now given by:
\begin{equation}
u \in \underset{|r| \leqslant u_{\max}}{\arg \max }\left(\beta_{2}.r^{2}+\beta_{1}.r\right)
\label{uoptreg}
\end{equation}
where:
\begin{align} \label{alfaoptregcoef}
&\beta_{2}:=-k \\
&\beta_{1}:=p_{\gamma}.\left(1+\lambda.\left(\dfrac{\rho(h)Sv}{2m}-1 \right) \right)
\end{align}
The transversality and Hamiltonian conditions remain unchanged:
\begin{equation}
\forall t \in [0, t_{f}], \hspace{0.2cm} H_{\lambda}^{k}(\xi(t),p(t),-1,u(t))=0
\label{Hamregtf}
\end{equation}
The transversality conditions on the adjoint state lead to:
\begin{equation}
p_{v}(t_{f})=0 \text{ and } p_{m}(t_{f})=0
\end{equation}
\subsection{Continuation on the dynamics and boundary conditions} 
The shooting function $S_{\lambda}: \mathbb{R}^{6} \mapsto \mathbb{R}^{6}$ is defined by:
\begin{equation}
S_{\lambda}: \left(\begin{array}{cc}
p(0) \\
t_{f}
\end{array}
\right) \longmapsto \left(\begin{array}{cccc}
x(t_{f})-x_{1}^{c}(\lambda) \\
h(t_{f})-h_{1}^{c}(\lambda) \\
\gamma(t_{f})-\gamma_{1}^{c}(\lambda) \\
p_{v}(t_{f}) \\
p_{m}(t_{f}) \\
H_{\lambda}^{k}\left(\xi(t_{f}),p(t_{f}),-1,u(t_{f})\right)
\end{array}
\right)
\label{S6}
\end{equation}
One performs the diagonal continuation on the following coordinates:
\begin{align*}
&h_{c}^{c}(\lambda)=h_{c}^{i}+\lambda(h_{c}-h_{c}^{i}) \\
&h_{0}^{c}(\lambda)=h_{c}^{c}(\lambda)+\lambda(h_{0}-h_{c}^{c}(\lambda)) \\
&v_{0}^{c}(\lambda)=v_{0}^{i}+\lambda(v_{0}-v_{0}^{i}) \\
&\gamma_{0}^{c}(\lambda)=\gamma_{0}^{i}+\lambda(\gamma_{0}-\gamma_{0}^{i}) \\
&x_{1}^{c}(\lambda)=x_{m}^{i}+\lambda^{2}(x_{f}-x_{m}^{i}) \\
&h_{1}^{c}(\lambda)=h_{c}^{i}(\lambda)+\lambda(h_{f}-h_{c}^{i}(\lambda)) \\
&\gamma_{1}^{c}(\lambda)=\gamma_{m}^{i}+\lambda(\gamma_{f}-\gamma_{m}^{i})
\end{align*}
where $\lambda \in [0, 1]$ is the homotopy factor increasing from $0$ to $1$. During this continuation, the maximum value of the control saturation was set to $u_{\max}^{i}=25$. The initalization values are the following: $h_{c}^{i}=0.5$, $x_{m}^{i}=5$, $\gamma_{0}^{i}=0$, $\gamma_{m}^{i}=0$, $v_{0}^{i}=1$. The remaining coordinates are unchanged with respect to their numerical values \eqref{boundcondbunt}. To ease the initialization of the shooting, we set initial and final values of $h$ at the cruise altitude: by proceeding so, the initial trajectory is likely horizontal at iso-cruise altitude $h_{c}$. The first guess $z(0):=(p_{x}(0), p_{h}(0), p_{v}(0), p_{\gamma}(0),p_{m}(0), t_{f})=(0.5,1.0,1.0,1.0,1.0, 6.0)$ enables a quick convergence.
\vspace{0.1cm}

\noindent The continuation step was set to the value $\delta \lambda=0.001$ and the regularization coefficient to $k_{\max}=100$. We highlight here that lower values of regularization parameter $k$ slow down significantly (or make fail) the algorithm. We infer that this observation is based on the following fact: when $k$ is ``too low'', then the central (level flight) arc becomes ``too close'' to a singular arc (that is along which the linearized system is not controllable, see definition \ref{singdef}). In this situation, it is well known that the continuation fails (see for instance Theorem 3 of \cite{TreZhuCer2017} for deeper expanations).  Consequently, $k$ can be considered as a ``singular perturbation'' in our optimal control problem. For deeper justification, please refer for instance to \cite{Ask2023}.
\vspace{0.1cm}

\noindent Under the previous settings, the processing time is around 71 seconds.  We display the optimal solution at the end of the first continuation:
\begin{figure}[h!]
\begin{flushleft}
\captionsetup{justification=centering}
\includegraphics[width=1.1\linewidth]{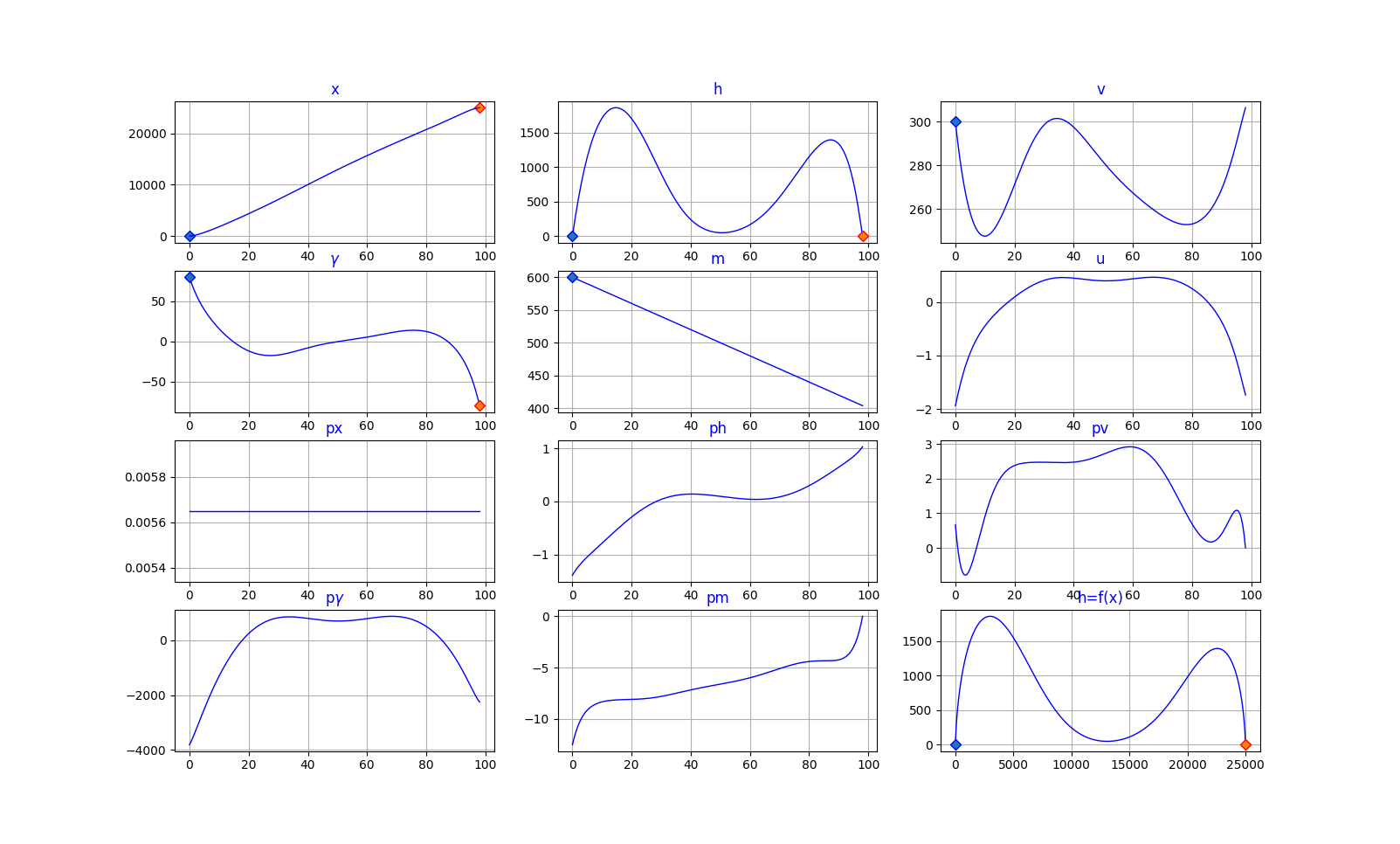}
\caption{Optimal state \& adjoint state for $\lambda=1$}
\label{fig7}	
\end{flushleft}
\end{figure}
\FloatBarrier
\noindent The following table provides the terminal accuracy:
\begin{table}[h!]
\begin{center}
\begin{tabular}{|c|c|}
 \hline
\textbf{Shooting output} & \textbf{Value} \\
\hline
$|x(t_{f})-x_{f}|$ & $7.10^{-10}$ \\
\hline
$|h(t_{f})-h_{f}|$ &  $2.10^{-10}$ \\
\hline
$|\gamma(t_{f})-\gamma_{f}|$ &  $5.10^{-11}$ \\
\hline
$p_{v}(t_{f})$ & $2.10^{-11}$ \\
\hline
$p_{m}(t_{f})$ & $2.10^{-11}$ \\
\hline
$H_{1}^{k}(t_{f})$ & $5.10^{-13}$ \\
\hline
\end{tabular}
\end{center}
\caption{Terminal accuracy}
\label{tab:shootstep1}
\end{table}

\section{From $(\textbf{OCP})_{k}$ to $(\textbf{OCP})$} \label{contphase2}
\subsection{Continuation strategy}
\noindent Once the boundary conditions reached, we implement a simultaneous decreasing continuation on the regularization coefficient $k$ and the control saturation $u_{\max}$. To do so, we implement a variant of the shooting method, which we will formally call "shooting from the middle" (refer to \cite{TreZua2015} and \cite{CaiFerTreZid2022} for deeper description and relevant examples). Using the notations of section \ref{indirect}, the shooting from the middle consists in choosing $z(t_{f}/2)$ as unknown and then:
\begin{itemize}
\item integrating backward in time the differential system over $[0, t_{f}/2]$ in order to get the value of $z(0)$
\item integrating forward in time the differential system over $[0, t_{f}/2]$ in order to get the value of $z(t_{f})$
\end{itemize}
This very simple variant appears to be very efficient for two main reasons:
\begin{itemize}
\item it divides the time horizon by a factor 2 by starting the shooting at time $t=t_{f}/2$
\item it is extremely useful when initialized at the turnpike (which is our case) as it exploits the local stability of the solution.
\end{itemize}
 Of course, "the price to pay" is the increase of the number of unknowns (and the conditions conditions to fullfill). 
 \vspace{0.2cm}

 \noindent We implement the shooting from the middle. The shooting is performed on the state $\xi\left(\frac{t_{f}}{2}\right)$ and costate  $p\left(\frac{t_{f}}{2}\right)$ taken at time $t_{f}/2$ and of course on the time of flight $t_{f}$. This continuation is parametrized by $\theta  \in [0, 1]$.
\subsection{Shooting function}
The shooting function $S_{\theta}: \mathbb{R}^{21} \mapsto \mathbb{R}^{21}$ is defined by:
\begin{equation}
S_{\theta}: \left(\begin{array}{cccc}
\xi^{-}\left(t_{f}/2\right) \\
\xi^{+}\left(t_{f}/2\right) \\
p^{-}\left(t_{f}/2\right) \\
p^{+}\left(t_{f}/2\right) \\
t_{f}
\end{array}
\right) \longmapsto \left(\begin{array}{cccc}
\xi(0)-\xi_{0} \\
x(t_{f})-x_{f} \\
h(t_{f})-h_{f} \\
\gamma(t_{f})-\gamma_{f} \\
p_{v}(t_{f}) \\
p_{m}(t_{f}) \\
H_{1}^{k^{c}(\theta)}\left(\xi(t_{f}),p(t_{f}),-1,u(t_{f})\right) \\
\xi^{-}\left(t_{f}/2\right)-\xi^{+}\left(t_{f}/2\right) \\
p^{-}\left(t_{f}/2\right)-p^{+}\left(t_{f}/2\right)
\end{array}
\right)
\label{S21}
\end{equation}
One performs the decreasing continuation over $k_{2}$ and $u_{\max}$:
\begin{align*}
&k^{c}(\theta)=k_{\max}+\theta(k_{\min}-k_{\max}) \\
&u_{\max}^{c}(\theta)=u_{\max}^{i}+\theta (u_{\max}-u_{\max}^{i})
\end{align*}
where $\theta \in [0, 1]$ is the homotopy factor increasing from $0$ to $1$. The initial guess is obtained by extracting $(\xi(t_{f}/2), p(t_{f}/2))$ and $t_{f}$ from the first step continuation for $\lambda=1.0$. After several numerical tests, we set $k_{\min}=2$ which is the minimum value for which  one obtains the convergence for $\theta=1.0$. We recall that $u_{\max}^{i}=25$  and $u_{\max}=2.0$ (nominal value). 
\vspace{0.1cm}

\noindent We set the continuation step $\delta \theta=0.005$ which ensures a reasonable compromise between the computation time ($150$ sec) and the algorithm accuracy. The solution obtained with the shooting method is presented hereafter and compared to the non regularized solution obtained with the direct method:
\begin{figure}[h!]
\begin{flushleft}
\captionsetup{justification=centering}
\includegraphics[width=1.\linewidth]{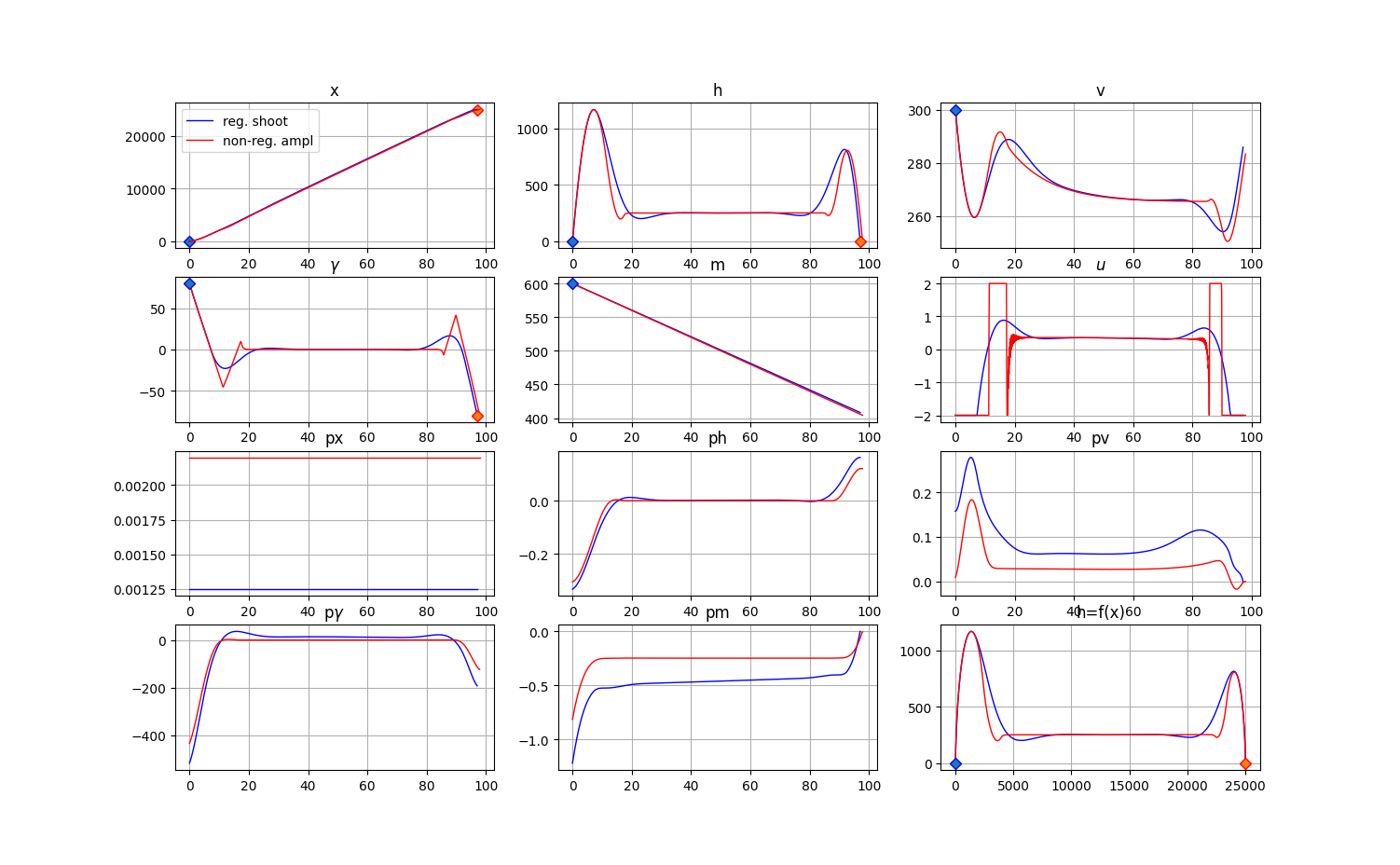}
\caption{Optimal state \& adjoint state for $\theta=1$}
\label{fig8}	
\end{flushleft}
\end{figure}
\FloatBarrier
The following table provides the accuracy of the continuation for $\theta=1.0$:
\begin{table}[h!]
\begin{center}
\begin{tabular}{|c|c|}
 \hline
\textbf{Shooting output} & \textbf{Value} \\
\hline
$|x(0)-x_{0}\|$ & $1.10^{-7}$ \\
\hline
$|h(0)-h_{0}|$ &  $7.10^{-7}$ \\
\hline
$|v(0)-v_{0}|$ &  $3.10^{-8}$ \\
\hline
$|\gamma(0)-\gamma_{0}|$ &  $2.10^{-8}$ \\
\hline
$|m(0)-m_{0}|$ &  $1.10^{-13}$ \\
\hline
$|x(t_{f})-x_{f}|$ & $3.10^{-6}$ \\
\hline
$|h(t_{f})-h_{f}|$ &  $5.10^{-7}$ \\
\hline
$|\gamma(t_{f})-\gamma_{m}|$ &  $2.10^{-7}$ \\
\hline
$p_{v}(t_{f})$ & $7.10^{-9}$ \\
\hline
$p_{m}(t_{f})$ & $4.10^{-7}$ \\
\hline
$H_{1}^{2}(t_{f})$ & $9.10^{-10}$ \\
\hline
$\|\xi^{-}\left(t_{f}/2\right)-\xi^{+}\left(t_{f}/2\right)\|$ & $0.0$\\
\hline
$\lVert p^{-}\left(t_{f}/2\right)-p^{+}\left(t_{f}/2\right) \rVert$ & $0.0$\\
\hline
\end{tabular}
\end{center}
\caption{Terminal accuracy}
\label{tab:shootstep1}
\end{table}
\FloatBarrier
\noindent The optimal cost of the non-regularized problem is around 210 whereas the corresponding cost (computed without the penalization term) of the regularized problem solution is around 230. Actually the difference is mainly due to the fact that the regularization induces an error in the tracking of the cruise altitude. However, from the operational point of view, the regularized solution is totally acceptable: indeed, the times of flight and the altitude overshoots at the beginning and the end are comparable.
\vspace{0.1cm}

\noindent Finally we highlight that, when initialized directly with the adequate guess (that is to say computed through the process of continuations detailed above), the shooting method converges quasi-instantaneously, providing the optimal solution.

\section{Conclusion}
\noindent The regularization of the initial control problem allows to successfully compute a solution. The key idea behind this is the increase of the regularization coefficient $k$ which enhances the well posedness of the shooting method and consequently ensures the convergence of the algorithm. The other advantage is that, contrary to the non-regularized optimization problem, one does not need to know the structure of the optimal trajectory prior to the implementation, which simplifies the analysis of the problem.
\vspace{0.1cm}

\noindent Moreover the regularized solution is relatively close to the original one in terms of cost and consequently totally acceptable from the operational point of view (see figure \ref{fig7}).
\vspace{0.1cm}

\noindent From the methodological point of view, we have exploited the (partial) turnpike phenomenon on the state variables $(h,\gamma)$ in order to implement efficiently a variant of the shooting method (initializing it ``at the middle''). As mentioned previously, this variant enhances the stability of the shooting as the control system remains ``invariant'' in the neighbourhood of the turnpike.
\vspace{0.1cm}

\noindent The computation time of the shooting method associated to the continuations can appear as relatively long. We highlight that, when initialized with the adequate guess, the shooting method provides the optimal solution \emph{quasi instantaneously}. As mentioned before, this is a major advantage of this approach with respect to the others. However, it requires to compute and store off-line a mapping of initial guesses corresponding to a grid of boundary conditions.
\vspace{0.1cm}

\noindent We have chosen to formulate the cruise altitude constraint through a penalization term in the cost. An alternative method is to introduce intermediate state constraints (for instance, $h(t)=h_{c}$ and $\gamma(t)=0$ for any $t \in [t_{1}, t_{2}]$ where $0<t_{1}<t_{2}<t_{f}$ are to be determined) and then compute an optimal solution using the hybrid Pontryagin maximum principle (see for instance \cite{DmiKag2008}-\cite{DmiKag2011}). However, we did not test this method in this paper.

\bibliographystyle{plain}
\bibliography{biblio}

\end{document}